\RequirePackage{ifpdf}
\ifpdf 
\documentclass[pdftex]{sigma}
\else
\documentclass{sigma}
\fi

\def\bbC{{\mathbb{C}}}
\def\bbR{{\mathbb{R}}}
\def\bbS{\mathbb{S}}
\def\bbV{{\mathbb{V}}}
\def\bbW{{\mathbb{W}}}

\def\cB{{\mathcal{B}}}
\def\cD{{\mathcal{D}}}
\def\cE{{\mathcal{E}}}
\def\cF{{\mathcal{F}}}
\def\cL{{\mathcal{L}}}
\def\cM{{\mathcal{M}}}
\def\cN{{\mathcal{N}}}
\def\cW{{\mathcal{W}}}
\def\cY{{\mathcal{Y}}}

\def\bul{{\scriptstyle\bullet}}
\def\ci{C^\infty}
\def\dc{\hbox{$\nabla\!\!\!\! /\;$}}
\newcommand{\minidcs}{\mbox{\boldmath{$\scriptstyle{{\dc}^2}$}}}
\newcommand{\minidcso}{\mbox{\boldmath{$\scriptstyle{{\dc}^2_0}$}}}
\def\dint{{\displaystyle{\int}}}
\def\dsum{{\displaystyle{\sum}}}
\def\endo{{\rm End}}
\def\eye{{\sf i}}
\def\kth{k^{\underline{\rm{th}}}}
\def\njs{|dJ|^2}
\def\prim{{\rm Prim}}
\def\sdint{{\dint\!\!\!\!\!{\rm-}}}
\def\tr{{\rm tr}}
\def\Tr{{\rm Tr}}
\def\dfrac{\frac}
\let\a=\alpha
\let\b=\beta
\let\d=\delta
\let\D=\Delta
\let\e=\varepsilon
\let\f=\varphi
\let\g=\gamma
\let\G=\Gamma
\let\h=\ell
\let\l=\lambda
\let\L=\Lambda
\let\m=\mu
\let\nd=\nabla
\let\N=\nabla
\let\p=\pi
\let\ptl=\partial
\let\s=\sigma
\let\u=\theta
\let\x=\xi
\let\w=\omega
\let\W=\Omega
\let\y=\psi
\let\z=\zeta

\newcommand{\nnn}[1]{(\ref{#1})}
\def\var#1{{#1}{}^\bul}

\begin{document}

\allowdisplaybreaks

\renewcommand{\PaperNumber}{090}

\FirstPageHeading

\renewcommand{\thefootnote}{$\star$}

\ShortArticleName{$Q$-Curvature, Spectral Invariants, and Representation Theory}

\ArticleName{$\boldsymbol{Q}$-Curvature, Spectral Invariants,\\ and Representation Theory\footnote{This paper is a
contribution to the Proceedings of the 2007 Midwest
Geometry Conference in honor of Thomas~P.\ Branson. The full collection is available at
\href{http://www.emis.de/journals/SIGMA/MGC2007.html}{http://www.emis.de/journals/SIGMA/MGC2007.html}}}

\Author{Thomas P.~BRANSON~$^\dag$}

\AuthorNameForHeading{T.P.~Branson}

\Address{$^\dag$~Deceased}

\URLaddress{\url{http://www.math.uiowa.edu/~branson/}}

\ArticleDates{Received August 01, 2007 from {\it Xingwang Xu}; Published online September 16, 2007}

\Abstract{We give an introductory account of functional determinants of
elliptic operators on manifolds and Polyakov-type formulas for
their inf\/initesimal and f\/inite conformal varia\-tions.
We relate this to extremal problems and to the $Q$-curvature on
even-dimensional conformal manifolds.
The exposition is self-contained, in the sense of giving references
suf\/f\/icient to allow the reader to work through all details.}


\Classification{58J52; 53A30}

\begin{flushright}
\begin{minipage}{15cm}\it
During May and June in 2004, the Institute of Mathematical Sciences, at the National University of Singapore, held a program on geometric partial dif\/ferential equations. The program had invited several instructional lecturers.  Professor Thomas Branson was one of them. Originally we had planned to publish the lecture notes of the instructional lecturers. This paper was submitted by Professor Branson for that purpose.  In fact this was the only paper we had received by the deadline. Hence we decided to cancel the plan for a~proceedings volume.  Just after that period, Thomas Branson unexpectedly passed away.  We held the paper without knowing what we could do with it.  When the editors of this proceedings volume invited me to submit an article, I realised that this would be an ideal place for Professor Thomas Branson's paper.  I immediately submitted the paper to editors of the current proceedings. I would like to take this opportunity to express my sincere appreciation to the editors for their help. \\[3mm]
\null \hfill                                   Xingwang Xu (National University of Singapore)\\
\null \hfill E-mail: \href{mailto:matxuxw@nus.edu.sg}{matxuxw@nus.edu.sg}
\end{minipage}

\end{flushright}

\section{The functional determinant}

In order to get a feel for the spectral theory of natural dif\/ferential
operators on compact manifolds, recall the idea of {\em Fourier
series}, where one attempts to expand complex functions on the unit
circle ${\bf S}^1$ in $\bbC$ in the form
\[
{\displaystyle\sum_{k=-\infty}^\infty} c_ke^{\eye k\u}
=c_0+{\displaystyle\sum_{j=1}^\infty}(a_j\cos j\u+b_j\sin j\u).
\]
The trigonometric series is an expansion in real {\em
eigenfunctions} of the {\em Laplacian} $\D=-d^2/d\u^2$ (the eigenvalue
being $j^2$).  The exponential series is an expansion in
eigenfunctions of the operator $-\eye d/d\u$, which is a square root
of the Laplacian; the eigenvalue is $k$.

Suppose we wish to f\/ind the {\em determinant} of the Laplacian on the
circle, or of its square root.  Taking the naive view that the
determinant should be the product of the eigenvalues in some spectral
resolution, we immediately have some problems -- we are really being
asked to take the product of the nonnegative integers.  We could
declare this to be 0 by convention (0 being one of the eigenvalues),
but the question of computing a
determinant like
\[
\det\left(-\dfrac{d^2}{d\u^2}+1\right)
\]
still persists -- is there any way to make sense of such determinants?

There are, of course, other characterizations of the determinant of a
f\/inite-dimensional linear transformation $T$, and we might try to
generalize one of these.  For example, there is the exterior algebra
def\/inition, which appears in abstract index notation (see below for
an explanation) as
\[
\det T=T^{a_1}{}_{[a_1}T^{a_2}{}_{a_2}\cdots T^{a_m}{}_{a_m]}\,.
\]
There is also a
{\em zeta function} def\/inition: let
\[
\z_T(s)=\Tr\,T^{-s},\qquad s\in\bbC.
\]
(The exponential functions $z\mapsto z^{-s}$ must be well-def\/ined on
spec$(T)$; for this, since spec$(T)$ is f\/inite, we just need
$0\notin{\rm spec}(T)$.)  The determinant is then
\begin{equation}\label{dz}
\det T=e^{-\z'_T(0)}.
\end{equation}
Indeed, \nnn{dz} follows from the calculation
\[
\left.\dfrac{d}{ds}\right|_{s=0}\Tr\,T^{-s}=-\Tr\,\log\,T.
\]
The particular branch cut used to compute $\log\,T$ af\/fects the value
of $\z'_T(0)$ (moving it by multiples of $2\p\eye$), but has no ef\/fect
on its exponential.

In fact, the zeta function approach is a fruitful one
for operators like the Laplacian on ${\bf S}^1$.  Before pursuing it
more deeply, however, let us look at still another approach.
If $T$ is symmetric and does not have 0 in its spectrum, then
\begin{equation}\label{detint}
\det T=\p^m\left(\int_{\bbR^m}e^{-(Tx,x)}dx\right)^{-2}.
\end{equation}
This is, in fact, the sort of formula that one tries to imitate
in quantum theory, with the so-called {\em functional integrals}
-- integrals over function space.  Since there are problems assigning
meaning to such integrals,
one tries to evaluate the inf\/inite-dimensional analogue of the
left side of \nnn{detint} instead.  Among the possible interpretations
(or {\em regularizations}) of the determinant, one may choose
the zeta function, and if the operator in question has certain properties,
this succeeds in assigning a value.

To see what zeta function regularization says about
$\det(\D)$
on ${\bf S}^1$, let us just brutally remove the 0 eigenvalue.
This makes the zeta function
\[
\z_\D(s)={\displaystyle\sum_{j=1}^\infty}2j^{-2s}=2\zeta_{{\rm Riem}}(2s),
\]
where $\z_{{\rm Riem}}$ is the Riemann zeta function.  As a result,
\[
\det\D=e^{-4\z'_{{\rm Riem}}(0)}=(2\pi)^2.
\]
If we could take the eigenvalue product interpretation of the determinant
literally, this would say that the product of the positive integers is
$\sqrt{2\p}$.

Let us recall the process which assigns meaning to the expression
$\z'_{{\rm Riem}}(0)$.
The
series $\sum\limits_{j=1}^\infty j^{-s}$ converges uniformly and absolutely on
the half-plane
for Re$(s)>1$.
The resulting holomorphic function has an analytic continuation
which may be constructed using Dirichlet series and the Mellin transform,
to a meromorphic function on $\bbC$, with a single simple pole at
$s=1$.

What we would like to explore here is a generalization of this thinking
to a situation in which:
\begin{itemize}
\itemsep=0pt
\item ${\bf S}^1$ is replaced by a compact $n$-dimensional smooth manifold;
\item the Laplacian is replaced by a dif\/ferential operator
with {\em positive definite leading symbol}.
\end{itemize}
Roughly speaking,
the same construction will go through in this setting,
and there will be a~{\em functional determinant}.

For example, consider the {\em Laplacian} of an $n$-dimensional
Riemannian manifold
$(M^n,g)$.  One has the {\em Riemannian connection} $\nd$, which
allows {\em covariant differentiation} of vector f\/ields, and in fact
of all kinds of tensors.  If $f$ is a smooth function, then
$\nd\nd f$ is a 2-tensor.  Contracting, one gets
\[
\D f:=-g^{ab}\nd_a\nd_bf=-\nd^a\nd_af
\]
in {\em abstract index notation}.  Here indices which are repeated,
once up and once down, denote contractions.  Indices are raised and
lowered using the metric, and the number and
position of indices just indicate a tensor's rank and covariance
type (rather than
any choice of local frame).  For more on abstract index notation, see
\cite{PR84}.

The {\em leading symbol} of a dif\/ferential operator $D$ is,
more or less, what one gets by writing~$D$ in local coordinates (and frames, if $D$ is vector bundle valued),
replacing each ${\partial}/{\partial{x^a}}$ by an~$\eye\x_a$,
and keeping only the terms with the maximal number of $\x$'s.
Here $\x$ is an indeterminate covector f\/ield, or one-form; the result
is some sort of function on the cotangent bundle of $M$.
A precise formula for the leading symbol may be given, without
making choices of local frames, as follows.
Suppose $D$ has order $k$, and
carries sections of a vector bundle $\bbV$ to sections of
a vector bundle~$\bbW$.
Then if $f$ is a smooth function and
$\f$ is a smooth section of $\bbV$,
\[
\s_k(D)(x,(df)_x)\f=
\displaystyle\lim_{t\to\infty}
\frac{(D(e^{\eye tf}\f))_x}{t^k},
\]
where the subscript $x$ denotes the value of a section (in this
case, of the cotangent bundle or $\bbW$) at $x$.
Since any covector $\x$ at $x$ is $(df)_x$ for some $f$, this
completely (and, as is easily checked,
consistently) determines $\s_k(D)(x,\x)\f$.

For example, in the case of the Laplacian on ordinary functions,
\[
\D f=-g^{ab}\ptl_a\ptl_b+{\rm(lower\ order)},
\]
so the leading symbol is
\[
\s_2(\D)(x,\x)=g^{ab}(x)\x_a\x_b\,.
\]
The ``2'' in the $\s_2$ notation just keeps track of
the order of the operator
(the maximal homogeneity in $\x$).
In this sense, the leading symbol of the Laplacian ``is'' the metric
(or, more precisely, the metric inverse $g^{ab}$).

Back in the general situation,
$\s_k(D)$ may be thought of as a section of the
vector bundle
\[
{\rm Hom}({\rm Symm}(T^*M\otimes\cdots\otimes T^*M)\otimes\bbV,\bbW),
\]
using the identif\/ication of $k$-homogeneous functions on a
vector space (here, the
cotangent
space~$T_x^*M$) with symmetric $k$-tensors on that space.
The construction of a {\em total} symbol is another matter.
One may view this as depending on choices (i.e., as a function
carrying atlases of $M$ which locally frame $\bbV$ and $\bbW$
to local total symbols) as in \cite{gilkeybook}.
Or, one can assume more geometric structure and
use it to
try to def\/ine total symbols invariantly, as in \cite{widom}.  This will not
be an explicit issue here, though it is an issue in the underlying
analysis of pseudo-dif\/ferential operators, which provides
the heat operators and complex powers of elliptic operators on which
zeta function regularization relies.
Roughly speaking, a total symbol will keep track of what happens to
the operator under local Fourier transforms.

The Laplacian's leading symbol has a precise positivity property: if
$\x\ne 0$, then
\[
g^{ab}\x_a\x_b>0,
\]
by the positive def\/initeness of the Riemannian metric.
We can speak of something like this in the vector bundle valued
case, say for a $\kth$-order dif\/ferential operator carrying sections
of a~vector bundle $\bbV$ to sections of the same bundle $\bbV$,
replacing positivity by
positive def\/initeness, at least if $\bbV$ has a positive def\/inite metric.
We say that $D$ has {\em positive definite leading symbol} if
for each $\x\ne 0$ in $T_x^*M$, for each $x\in M$, we have
\[
\s_k(D)(x,\x)\ \ \mbox{{\rm positive\ def\/inite\ in}}\ \ \endo(\bbV_x).
\]
A positive def\/inite endomorphism $A$ is, of course, one with
$(Av,v)>0$ for each nonzero vector~$v$, where $(\cdot,\cdot)$ is
the inner product (provided here by the metric on the bundle $\bbV$).

Note that only even-order operators can have positive def\/inite leading
symbol, since if $k$ is odd and $\s_k(D)(x,\x)$ is positive def\/inite,
then $\s_k(D)(x,-\x)$ is negative def\/inite.

Perhaps the easiest examples of operators with positive def\/inite leading
symbol which do {\em not} have the form
\[
\D^h+{\rm(lower\ order)}
\]
are the dif\/ferential form operators $\d d+ad\d$, where $d$ and $\d$
are the exterior derivative and coderivative, and $a$ is a positive
constant not equal to 1.
To be more specif\/ic, if $\f$ is a $p$-form, then $d\f$
is the $(p+1)$-form
\begin{equation}\label{defd}
(d\f)_{a_0\ldots a_p}=(p+1)\nd_{[a_0}
\f_{a_0\ldots\hat{a}_s\ldots a_p]},
\end{equation}
and $\d\f$ is the $(p-1)$-form
\[
(\d\f)_{a_2\ldots a_p}=-\nd^b\f_{ba_2\ldots a_p}.
\]
In fact, one does not need a connection to def\/ine $d$; the right side
of \nnn{defd} is invariant as $\nd$ ranges over all symmetric af\/f\/ine
connections.
$\d$ does depend on the connection $\nd$, and in fact
is the formal adjoint of $d$ in the form metrics
\[
(\f,\y)_p=\dfrac1{p!}\f^{a_1\ldots a_p}\y_{a_1\ldots a_p}.
\]
This means that if $\f$ is a $p$-form and $\eta$ is a $(p+1)$-form
on a Riemannian manifold $(M,g)$, then
\[
\int_M(d\f,\eta)_{p+1}dv_g=\int_M(\f,\d\eta)_pdv_g
\]
provided either $\f$ or $\eta$ has compact support.
Here $dv_g$ is the Riemannian measure.  The formal adjoint property
makes it clear that $\d$ (unlike $d$) will generally change if we
vary the metric $g$.

The {\em form Laplacian}
\[
\D_{{\rm form}}=\d d+d\d
\]
is related to the Laplacian $\N^*\N$ (sometimes called the {\em
Bochner Laplacian} in this context) by the {\em Weitzenb\"ock formula}
\[
\d d+d\d=\N^*\N+\cW.
\]
Here $\cW$, the {\em Weitzenb\"ock operator}, is an action of the
curvature tensor which has order 0 as a~dif\/ferential operator.
As a result,
\[
\s_2(\D_{{\rm form}})(x,\x)=|\x|^2{\rm Id}_{{\rm form}}\,,
\]
where the length is measured by the metric: $|\x|^2=g^{ab}\x_a\x_b$.

But if $D$ has order $k$ and $E$ has order $\h$, then
\[
\s_{k+\h}(DE)(x,\x)=\s_k(D)(x,\x)\s_\h(E)(x,\x).
\]
In addition, formal adjoints (denoted by
the superscript $*$) exist and are locally determined, and
\[
\s_k(D)(x,\x)^*=\s_k(D^*)(x,\x).
\]
Thus $\d d+ad\d$ has positive semidef\/inite leading symbol for
$a\ge 0$, and for $a>0$, the calculation
\[
\left(\d d+ad\d\right)\left(\d d+a^{-1}d\d\right)=\D_{{\rm form}}^2
\]
(based on $dd=0$) shows that $\d d+ad\d$ has invertible leading symbol.
Thus $\d d+ad\d$ has positive def\/inite leading symbol for $a>0$.

Beyond dif\/ferential forms, in other natural bundles, the idea of
positive def\/inite leading symbol makes sense, but the calculus is
not as simple -- in particular, there is usually no complete analogue of the
relation $dd=0$.

To construct the zeta function and the {\em heat operator} of a
dif\/ferential operator $D$ on sections of $\bbV$, we make the following
assumptions.

{\bf Analytic assumptions.} $(M^n,g)$ is an $n$-dimensional
compact smooth manifold, and $\bbV$ is a~vector bundle over
$M$ with a Hermitian metric $h$.
$A$ is a smooth-coef\/f\/icient
dif\/ferential operator of positive order
on sections $\bbV$
which is {\em formally self-adjoint}:
\[
\int_Mh(\f,A\y)\,dv_g=\int_Mh(A\f,\y)\,dv_g,
\]
and $A$ has positive def\/inite leading symbol.

From these assumptions we may conclude that
$A$ has even order $2\ell$, since taking $\x\mapsto-\x$
reverses the sign of the leading symbol for odd-order operators.
In addition, $L^2(M,\bbV)$ has a~complete
orthonormal set $\{\f_j\}$ of eigensections for $A$
with real eigenvalues
\[
\l_0\le\l_1\le\cdots\le\l_j\le\cdots\,,
\]
satisfying
\begin{equation}\label{weylasym}
\l_j\sim{\rm const}\cdot j^{2\ell/n},\qquad j\to\infty.
\end{equation}
\nnn{weylasym} is the {\em Weyl asymptotic} law; see, e.g., \cite{gilkeybook}.

If $r\ge 0$, the {\em Sobolev space} $L^2_r(M,\bbV)$
is the domain of the operator
\[
\displaystyle
(A-\l_0+1)^{r/2\ell}:\quad \sum_jc_j\f_j\mapsto\sum_jc_j(\l_j-\l_0+1)^{r/2\ell}\f_j
\]
in $L^2(M,\bbV)$.  The absence of any reference to $A$ in the
notation $L^2_r(M,\bbV)$ suggests that this domain is independent of
$A$, and in fact it is.  (As a consequence of the compactness of
$M$, the Sobolev spaces do not even depend on the metric of
the tangent bundle or of
$\bbV$; see \cite{palais}.)
In addition, there is always some usable operator $A$ -- for example,
make metric and connection choices and
take the
Laplacian $\N^*\N$.
Though the Sobolev space $L^2_r(M,\bbV)$ is
Hilbertizable, it is more productive to think of it as a Banach space,
with the equivalence class of norms represented by
\[
\|\f\|^2=\int_Mh((A-\l_0+1)^{r/2\ell}\f,(A-\l_0+1)^{r/2\ell}\f)dv_g.
\]
This is just a more-or-less random choice of a norm from the
relevant equivalence class.  When delicate invariance considerations
are in play, it sometimes pays to pick an invariant representative
of this class; see, for example, the last section below.

The dual space is denoted $L^2_{-r}(M,\bbV)$; for real numbers $r\le s$,
\begin{gather*}
({\rm distributional\ sections})=
\cE(M,\bbV)\subset L^2_r(M,\bbV)\subset L^2_s(M,\bbV)\\
\phantom{({\rm distributional\ sections})=}{}
\subset
\ci(M,\bbV)=({\rm smooth\ sections}),
\end{gather*}
and $L^2=L^2_0$.
The Sobolev Lemma shows that
\[
\displaystyle\ci=\displaystyle\bigcap_rL^2_r,
\]
and this, together with the estimates
def\/ining the space of distributions, shows that
\[
\displaystyle\cE=\displaystyle\bigcup_rL^2_r.
\]

The {\em heat operator} $e^{-tA}$, $t>0$, is of trace class, and is
{\em smoothing}: it carries $L^2_r$ continuously to $L^2_s$ for any
$r,s$, and thus carries distributional sections to to smooth sections.
As a result,
composition with $e^{-tA}$ has a very ``civilizing'' ef\/fect on
other operators: if, for some $N$, the operator $B$ carries $L^2_r$ to
$L^2_{r-N}$ for each $r$, then $Be^{-tA}$ and $e^{-tA}B$ carry
distributional sections to smooth sections.

The {\em heat trace} has {\em small-time asymptotics}:
\[
\Tr_{L^2}e^{-tA}\sim
a_0[A]t_{}^{\frac{-n}{2\ell}}+a_1[A]t^{\frac{-n+1}{2\ell}}
+\cdots+
a_k[A]t^{\frac{-n+k}{2\ell}}+\cdots,\qquad t\downarrow 0.
\]
This means that if we take the sum on the right out
to the $\kth$ term, we have accuracy to
the order of the next term: $O(t^{(-n+k+1)/2\ell})$.

Note that one can conclude this accuracy immediately from a much
weaker statement -- accuracy to order $O(t^{N(k)})$, where
$N(k)\to\infty$ as $k\to\infty$.  To see this, pick $k_1$ so
that $N(k_1)\ge(-n+k+1)/2\ell$, and expand out to order $k_1$.
Then the error in approximating
$\Tr_{L^2}e^{-tA}$
by
\[
a_0t_{}^{\frac{-n}{2\ell}}+a_1t^{\frac{-n+1}{2\ell}}
+\cdots+
a_kt^{\frac{-n+k}{2\ell}}
\]
is
\[
a_{k+1}t^{(-n+k+1)/2\ell}+\cdots+
a_{k_1}t^{(-n+k_1)/2\ell}+O\big(t^{(-n+k+1)/2\ell}\big)
=O\big(t^{(-n+k+1)/2\ell}\big).
\]

The way in which the $a_k$ depend on (the total symbol of)
$A$ can be minutely analyzed, and it is useful to do so.
One qualitative observation is that the
$a_{{\rm odd}}$ vanish in the current setting.
(On manifolds with boundary, they consist entirely of boundary
integrals.)
Again, an excellent reference is \cite{gilkeybook}.

In addition to the asymptotic series above, the heat operator
trace also expands as a
{\em Dirichlet series of the
first kind},
\[
\Tr_{L^2}(e^{-tA})=\sum_{j=0}^\infty e^{-\l_jt}\qquad(t>0).
\]
For example, when $A$ is the operator $-d^2/d\u^2$
on the circle, the series is
\[
1+\sum_{j=1}^\infty 2e^{-j^2t}.
\]

The corresponding Dirichlet series of the {\em second kind} take
the form
\[
\sum_{j=0}^\infty\l_j^{-s}.
\]
Actually, this is not entirely accurate, since
\nnn{weylasym} allows a f\/inite number of $\l_j$ to be negative, and
a f\/inite number to be 0.  Some sort of artif\/icial ``f\/ix'' is thus needed
to write such a Dirichlet series; we shall make the choice
\[
\displaystyle\sum_{\l_j\ne 0}|\l_j|^{-s};
\]
this is the {\em zeta function} $\z_A(s)$.
Whatever artif\/icial convention we choose for dealing with nonpositive
eigenvalues,
the ef\/fects will be felt in all succeeding formulas
which make contact with the zeta function.

The two types of Dirichlet series are related by the {\em
Mellin transform}
\[
(\cM F)(s)=\dfrac1{\G(s)}\int_0^\infty t^{s-1}F(t)dt,
\]
under which
\begin{equation}\label{mel}
\cM:e^{-\l t}\mapsto\l^{-s}\qquad(\l>0).
\end{equation}
Because of \nnn{mel},
\[
\cM:\underbrace{\Tr\,e^{-tA}}_{=:Z_A(t)}\mapsto
\underbrace{\Tr\,A^{-s}}_{\z_A(s)}
\]
when $A$ has positive spectrum.
When there are nonpositive eigenvalues, the convention on
$\z_A(s)$ adopted above forces us to replace $Z_A(t)$ by
\begin{equation}\label{fixx}
\tilde Z_A(t)=\sum_{\l_j\ne 0}e^{-t|\l_j|}=
Z_A(t)-q[A]+2\displaystyle\sum_{\l_j<0}\sinh(\l_jt),
\end{equation}
where $q[A]$ is the multiplicity of 0 as an eigenvalue of $A$.

The analytic continuation of the zeta function, which so far is only
well-def\/ined on ${\rm Re}\,s>n/2\ell$ (by the Weyl asymptotics \nnn{weylasym}),
is accomplished by substituting the small-time heat asymptotics into
the Mellin transform expression for the zeta function.
More precisely, assuming for the moment that there are
no nonpositive eigenvalues,
\[
\G(s)\z_A(s)=\sum_{k=0}^Na_k[A]
\left(s-\dfrac{n-k}{2\ell}\right)^{-1}
+\underbrace{\dint_0^1t^{s-1}O(t^{(-n+k+1)/2\ell})dt}_{{\rm regular,\
Re}\,s>\frac{n-k}{2\ell}}
+\underbrace{\dint_1^\infty t^{s-1}\sum_je^{-t\l_j}dt}_{\rm entire}.
\]
Since we can do this for arbitrarily large $N$, the $\G\z$ function
has been analytically continued to a meromorphic function with
possible simple poles at the $(n-k)/2\ell$, where $k$ is an {\em even}
natural number.
A useful viewpoint is that we are ``approximating'' the
$\G\z$ function; the error is ``small'' when it
is regular on a half plane starting far to the left.
In this sense the ``smallest'' functions are the entire ones.

To paraphrase the above, the $\G\z$ function has at most simple
poles on the real axis,
at equal decrements starting with $s=n/2\ell$.  Thus the zeta function
has no poles at nonpositive integers~-- they're resolved by
the zeros of $1/\G(s)$.  In particular,
\[
\z_A(s)\ {\rm is\ regular\ at}\ s=0.
\]
This is what we needed in order to
state the zeta function def\/inition of the determinant,
\[
\det\,A=e^{-\z'_A(0)}.
\]

\begin{remark} The regularity of the zeta function at $s=0$
is a more delicate property than
the Weyl asymptotics \nnn{weylasym}.
It is possible to construct sequences of $\l_j$ satisfying the
Weyl asymptotics for which the poles of the corresponding zeta
functions perform various tricks -- in particular, one can place
a pole at $s=0$.
Somehow, regularity at $s=0$ indicates
that an ``organic'' property of dif\/ferential (or
pseudo-dif\/ferential) operators is
being detected by the zeta function.
\end{remark}

\begin{remark}
If we don't necessarily have positive eigenvalues,
the new spectral function $\tilde Z_A(t)$ of \nnn{fixx}
has its own modif\/ied heat coef\/f\/icients $\tilde a_k$,
def\/ined by
\[
-q[A]+2\dsum_{\l_j<0}
\sinh(t\l_j)+\sum_{k=0}^\infty a_k[A]t^{(k-n)/2}
\sim\dsum_{k=0}^\infty\tilde a_k[A]t^{(k-n)/2}.
\]
Note that $\tilde a_k=a_k$ for $k<n$.
To build more determinant-like properties
into the functional determinant, we might insert a sign to better
monitor the negative eigenvalues:
\[
\det A:=(-1)^{\#\{\l_j<0\}}e^{-\z'_A(0)}.
\]
\end{remark}

\section{Variation of the functional determinant}

We now begin to
imagine wiggling the operator $A$ in various ways.
For example, we could let the operator $A$ depend on the {\em geometry}
(i.e.\ the metric $g$)
like the Laplacian does, and vary the metric.
An important way to vary the metric is {\em conformally}:
\[
\bar g=\W^2g,\qquad 0<\W\in\ci(M).
\]
In fact, by writing
\[
g_\w=e^{2\w}g_0,
\]
we can parameterize the {\em conformal class} of $g_0$ by $\w\in\ci(M)$.

There are dif\/ferential operators that are well adapted to this variation,
namely
the {\em conformal covariants}, or {\em conformally invariant differential
operators}.  Such an ``operator'' is really a rule assigning
operators to metrics in a natural way, and has
\begin{equation}\label{cco}
\bar g=\W^2g\ \Rightarrow\ \bar D\f=\W^{-b}D(\W^a\f)
\end{equation}
for some {\em conformal biweight} $(a,b)$ (and all
sections $\f$ of the appropriate bundle).
By rephrasing things in terms of operators on
{\em density bundles} of appropriate weights, we could restate
this as $\overline{D'}=D'$, for a uniquely determined operator $D'$,
but we shall not pursue that viewpoint here.
The best-known example of a conformal covariant is the {\em conformal
Laplacian}
\[
Y=\D+\dfrac{n-2}{4(n-1)}K\qquad (K={\rm scalar\ curvature}),
\]
on scalar functions; its conformal biweight is $((n-2)/2,(n+2)/2)$.

The {\em infinitesimal form} of the conformal covariance
relation is equivalent to the f\/inite form~\nnn{cco}, and
is sometimes more useful.  Let the metric run through the curve
$\{g_{\e\w}\}_{\e\in\bbR}$
within a conformal class, and let
\[
\var{\ }\;=\dfrac{d}{d\e}\Big|_{\e=0}.
\]

Then a conformal covariant has
\[
\var{D}\f=aD(\w\f)-b\w D\f,
\]
or, in operator terms,
\[
\var{D}=-(b-a)\w D+a[D,\w].
\]
(The f\/inal $\w$ in this formula should be interpreted as multiplication
by $\w$.)
Inf\/initesimal conformal covariance is equivalent to f\/inite conformal
covariance, since any metric $g_\w$
in the conformal class may be connected to $g_0$ by the
curve $\{g_{\e\w}\}_{\e\in[0,1]}\,$.

In addition to the metric
$g$, there are other standard ingredients involved in some geometric
operators -- the volume form $E$, and/or the
fundamental tensor-spinor $\g$.
If these are used, they should be scaled compatibly in forming $\bar D$:
\[
\bar E=\W^nE,\qquad\bar\g=\W^{-1}\g.
\]
An example of a conformal covariant for which this is relevant
is the {\em Dirac operator} $\dc=\g^a\N_a$ on spinors.
its conformal biweight is $((n-1)/2,(n+1)/2)$.
Another example is the operator $\star d$ on $(n-1)/2$-forms in
odd dimensions; here the {\em Hodge star} operator carries
$p$-forms to $(n-p)$-forms on an oriented manifold, and is def\/ined by
\[
(\star\f)_{a_{p+1}\ldots a_n}=\dfrac1{p!}
E_{a_1\ldots a_n}\f^{a_1\ldots a_p}.
\]
This makes $E$ enter the def\/inition of $\star d$ on $\L^{(n-1)/2}$; the
conformal biweight of this operator is $(0,2)$.

We would like to know what happens to our spectral data if we vary
conformally, and if our operator $A$ depends on $g$ in a conformally
``nice'' way.  It is generally hopeless to try to track the motion of
individual eigenvalues (except in very special circumstances; see the
discussion after~\nnn{otherinvts}
below).  However, it is sometimes possible to track the motion of
quantities like the heat coef\/f\/icients $a_k$ and the functional
determinant.  The reason is that these quantities have natural
def\/initions in terms of traces of operators; such def\/initions take
account of the whole spectrum.  If the variations of the operators
are known, there is some chance of computing the variation of the
natural spectral quantity.  As we shall see, conformal variation of
suitable operators yields a setting in which variational computations
can be done.

The good class of suitable operators is given by:

\medskip

\noindent
{\bf Conformal assumptions.}
$A=D^h$ is a positive power of a positive-order conformal cova\-riant~$D$.

\medskip

For example, the conformal Laplacian satisf\/ies both the analytic
and conformal assumptions, since it is itself a conformal covariant.
The Dirac operator $\dc$ is conformally covariant, but
does not have positive def\/inite leading symbol.  However its
{\em power} $\dc^2$, though not conformally covariant,
has positive def\/inite leading symbol (in fact,
leading symbol $|\x|^2\,{\rm Id}$), and so satisf\/ies our analytic and
conformal assumptions.  The operator $\star d$ on $\L^{(n-1)/2}$
does not have invertible leading symbol, so there is no hope of
any of its powers satisfying the analytic assumptions.

Under our analytic and conformal assumptions, we have the following:

\begin{theorem}[Theorem on Variations \cite{bo86,bo91,pr87,dk78}]
With notation as explained just below,
\begin{gather}
\var{a_k[A]}=(n-k)a_k[\w,A]=(n-k)\int_M\w U_k[A]dv_{g_0},\label{cit} \\
\var{\z'_A(0)}=2\ell\int_M\w\left\{U_n[A]-\sum_{\l_j=0}|\f_j|^2
\right\}dv_{g_0},\label{vardet}
\end{gather}
where $\{\f_j\}$ is any orthonormal basis of the $0$-eigenspace of $A$.
\end{theorem}

The \cite{pr87} and \cite{dk78} references really just treat the
conformal Laplacian.

The statement uses the {\em local} heat coef\/f\/icients $a_k[\w,A]$
and $U_k[A]$.  The f\/irst of these, $a_k[\w,A]$ is a term in the
small-time asymptotics of $\Tr_{L^2}\w e^{-tA}$:
\[
\Tr_{L^2}\w e^{-tA}\sim a_0[\w,A]t_{}^{\frac{-n}{2\ell}}+a_1[\w,A]
t^{\frac{-n+1}{2\ell}}
+\cdots+
a_k[\w,A]t^{\frac{-n+k}{2\ell}}+\cdots,\qquad t\downarrow 0.
\]
An analysis of the dependence of the heat coef\/f\/icients $a_k[\w,A]$
on $A$ and $\w$ shows that they are integrals of
dif\/ferential polynomials in
the total symbol of $A$, with coef\/f\/icients that depend (in a~universal way) smoothly
on the leading symbol and linearly on the jets of $\w$.  As a result,
after integration by parts,
\[
a_k[\w,A]=\int_M\w U_k[A]dv_g
\]
for some functions $U_k[A]$ built in a universal way from the total
symbol of $A$.

In particular, we have:

\begin{theorem}[Conformal Index Theorem \cite{bo86}]
Under the above assumptions,
the coefficient~$a_n[A]$ is conformally invariant.
\end{theorem}

Some other conformally invariant
quantities which are important in these calculations are
\begin{gather}
q[A]=\dim\cN(A)=\#\{\l_j=0\}, \qquad
\z_A(0)=a_n[A]-q[A], \qquad
\#\{\l_j< 0\}.\label{otherinvts}
\end{gather}
The invariance of $q[A]$ is immediate from the conformal covariance
relation.  The invariance of $\#\{\l_j< 0\}$ then follows from
this and Browder's Theorem \cite{browder}.  At the metric $g_0$,
for some $\d>0$, the negative eigenvalues are contained in
the interval $(-\infty,-\d)$, and the positive ones
in $(\d,\infty)$.
With the invariance of the number of zero eigenvalues, this is good
enough to keep the number of negative eigenvalues constant on
$\{g_{\e\w}\}$, for f\/ixed $\w$ and $|\e|$ bounded by some $\e_0>0$.
But this in turn implies constancy on the whole conformal class.

The local data $a_k[\w,A]$, $U_k[A]$ carry more information
than $a_k[A]$; in fact
\[
a_k[A]=a_k[1,A]=\int U_k.
\]
This more ref\/ined data (involving $\w$) is also {\em spectral},
but the spectra of many more operators (than just $A$)
are now involved.

If $A$ is {\em natural}, that is, if it is
built in a universal way from the metric (and possibly $E$ and/or~$\g$),
then the total symbol of $A$ must also be.
But then, in turn, the
$U_k[A]$ are also built in this way.
By {\em Weyl's invariant theory} \cite{weyl39}, the $U_k[A]$ are
built polynomially from the Riemann curvature and its
iterated covariant
derivatives.

For example, if $K$ is the scalar curvature and $\D$ is the Laplacian
on functions,
\begin{gather}
U_0[\D+aK]=(4\pi)^{-n/2},\nonumber \\
U_2[\D+aK]=(4\pi)^{-n/2}\left(\tfrac16-a\right)K, \label{DplusaK}\\
U_4[\D+aK]=(4\pi)^{-n/2}\cdot\tfrac1{180}\cdot
\left[90\left(\tfrac16-a\right)^2K^2
-|r|^2+|R|^2-30\left(\tfrac15-a\right)\D K\right],\nonumber
\end{gather}
where $R$ (resp.\ $r$) is the Riemann (resp.\ Ricci) tensor.
One reason that formulas like this are important is
hinted at by \nnn{vardet}: we might be able to parlay
explicit knowledge of $U_n[A]$ into information about
$\det\,A$ as a function on the conformal class $\{g_\w=e^{2\w}g_0\}$.
In fact, there are two genuinely dif\/ferent ways of doing this, one
based on \nnn{vardet}, and one based on both \nnn{vardet} and \nnn{cit}.

The f\/irst method is to simply
integrate the variation along the
curve $\{g_{\e\w}\}$, starting from $g_0$ and ending at $g_\w$.
As it turns out, the homogeneities are such that this is just
the integration of a polynomial in $\e$.  Specif\/ically,
\begin{equation}\label{cvar}
U_k[A_\w]=e^{-2k\w}\left(U_k[A_0]+U_k^{(1)}[A_0](\w)+\cdots
+U_k^{(k)}[A_0](\w)\right),
\end{equation}
where $U_k^{(i)}$ is a polynomial in the Riemann curvature and $\w$
and their covariant derivatives, and is of homogeneous degree $i$ in
$\w$.  (In fact, it depends only on $d\w$; see \cite{tbms85}.)
Here we have adopted the convention of labelling geometric
objects depending on the metric by the subscript $\w$, rather than
the subscript $g_\w$; for example,
\[
A_\w=A_{g_\w}.
\]
On the other hand,
\[
dv_\w=\e^{n\w}dv_0.
\]
As a result, $U_n[A_{\e\w}]dv_{\e\w}$ is polynomial in $\e$, of degree
at most $n$.

(By the conformal invariance of the conformal index
$\int U_n[A]dv$, we may
actually say ``of degree at most $n-1$''.  If $n$ is odd, $U_n[A]$
is identically zero, while if $n$ is even, homogeneity considerations
force the $U_n^{(n)}[A]$ term in $\nnn{cvar}$ to take the form
$c\cdot|d\w|^n$, for some universal constant $c$.  But if $c$
were nonzero, a simple scaling argument in $\w$ shows that
the conformal index is unbounded on the conformal class, contradicting
the fact that it is constant.)

The second method is an adaptation of the
physicists' {\em dimensional regularization}.
For simp\-licity in this discussion, suppose that $A$ has no
zero modes (eigenvalues) at $g_0$ (and thus on the whole conformal class).
Suppose that we are interested in the functional determinant in
dimension $n$.
The idea is to get a formula for $a_n[A]$ which is valid in an
inf\/inite number of dimensions $n'$, including $n'=n$.
Denote by Prim (for ``primitive'') the conformal anti-variation
corresponding to the conformal variation $\var{\ }$ above;
this is well-def\/ined up to a ``constant of integration'' -- a constant
function on the conformal class.
The Conformal Index Theorem then says that
\[
\prim\dint\w U_n[A]=\left\{
\begin{array}{ll}
\tfrac{a_n[A]}{n'-n}+{\rm const},\qquad & n'\ne n,\vspace{1mm} \\
\tfrac1{2\ell}\z'_A(0)+{\rm const},\qquad & n'=n.
\end{array}\right.
\]
The idea is now that the f\/irst formula continues to the second.
To standardize the constants of integration, we view the primitives
as functions of an indeterminate conformal factor $\w$, measured
from a base metric $g_0$, and require that our primitives vanish
at $g_0$:
\[
\prim_0\dint\w U_n[A]=\left\{
\begin{array}{ll}
\tfrac{a_n[A_\w]-a_n[A_0]}{n'-n},\qquad & n'\ne n, \vspace{1mm}\\
\tfrac1{2\ell}(\z'_{A_\w}(0)-\z_{A_0}(0)),\qquad & n'=n.
\end{array}\right.
\]
The analytic continuation to
the special
dimension $n'=n$ is then expressed by the formula
\[
\dfrac{a_n[A_\w]-a_n[A_0]}{n'-n}\Bigg|_{n'=n}=
\dfrac1{2\ell}(\z'_{A_\w}(0)-\z'_{A_0}(0)).
\]
As a result, we only need to know the {\em integrated} invariant
$a_n[A]$, not the whole local inva\-riant~$U_n[A]$.  But we need
to know it in ``all'' dimensions, or at least an inf\/inite number
of dimensions.  If we have this information, we get a formula for
the {\em difference} of log-determinants at two conformal metrics;
that is, the {\em quotient} of determinants.

Note that the f\/irst method (f\/ixing the dimension and
integrating along a curve in the conformal class)
is also giving us determinant {\em quotients}, since
the def\/inite integral described above computes
$\z'_{A_\w}(0)-\z'_{A_0}(0)$.

Of course, one cannot continue a {\em formula} from an inf\/inite number
of values of the independent variable to a special value.  But such a
continuation {\em is} possible in the case of a {\em rational}
formula.  That is, a rational function on $\bbC$ may be continued from
any inf\/inite set of values, because a polynomial function can.
Alternatively, a meromorphic function on $\bbC$ for which $\infty$ is
not an essential singularity may be continued from any set having
$\infty$ (or any other point of the Riemann sphere, for that matter)
as a limit point.  What we need is a way of phrasing the above
statements in terms of linear combinations, with rational
coef\/f\/icients, of a stable (under change of dimension) basis of local
invariants.  The ingredients out of which these invariants will be
built will sometimes include the indeterminate conformal factor $\w$
-- like the metric or the Riemann curvature, it will be thought of as
just another formal variable, rather than a specif\/ic function on a
specif\/ic manifold.

Of course one knows that the heat invariants are usually {\em not}
rational in the dimension.  For example, for the operators $\D+aK$ on
scalar functions of \nnn{DplusaK}, it looks as if $U_k[A]$ will have
the form $(4\p)^{-n/2}V_k[A]$ with $V_k[A]$ rational in $n$ (as long
as the constant $a=a(n)$ is rational in~$n$).  This is in fact true
for all $k$, not just the $k\le 4$ cases in \nnn{DplusaK}.  But the
exponential function~$(4\p)^{-n/2}$ certainly has an essential
singularity at $n=\infty$.  The remedy for this is the fact that the
$(4\p)^{-n/2}$ factor is {\em universal} -- one can factor it out of
the coef\/f\/icient of each invariant, and it's constant with respect to
the relevant variation of the geometry.  Thus it's really the
coef\/f\/icients of $V_k[A]$ (and related quantities) that we're talking
about, and these coef\/f\/icients are rational.  The possibility of
factoring out a universal function of $n$ so that we are left with
rational functions of $n$ is really a question about the form of the
dependence of the operator's leading symbol on the metric.  Similar
results hold for other leading symbols; for example~$|\x|^4$~\cite{tbcmp96}, and bundle-valued symbols which are not just powers of $|\x|^2$ times the identity~\cite{ab1}.

{\bf  Demonstration of the Theorem on Variations.}
The key to the variational formulas is the
conformal variation of the heat operator trace, or of the zeta function:
\begin{gather*}
\var{Z(t)}=
(\Tr\,e^{-tA})^\bul
 =-t\,\Tr\,\var{A}e^{-tA} =-t\,\Tr\,\underbrace{\var{D}}_{aD\w-b\w D}
hD^{h-1}\exp(-tD^h) \\
\phantom{\var{Z(t)}}{}=th(b-a)\,\Tr\,\w D^h\exp(-tD^h).
\end{gather*}
Note that the step $(\Tr\,e^{-tA})^\bul
=-t\,\Tr\,\var{A}e^{-tA}$ is not valid on the operator level, even
formally; the operators involved only have the same trace.  There is
no problem with domains due to the ``civilizing'' ef\/fect of
the smoothing operator $\exp(-tA)$~-- placing it next to a
pseudo-dif\/ferential operator
of f\/inite order produces an
operator of order $-\infty$, and {\em a fortiori} one that is bounded
on $L^2$.
By homogeneity considerations,
\[
h(b-a)=2\ell,
\]
so
\[
(\Tr\,e^{-tA})^\bul=2\ell t\,\Tr\,\w Ae^{-tA}
=-2\ell t\dfrac{d}{dt}\,\Tr\,\w e^{-tA}.
\]

However,
\[
t\dfrac{d}{dt}\big(t^{(k-n)/2\ell}\big)=\dfrac{k-n}{2\ell}
t^{(k-n)/2\ell}.
\]
So, provided that term-by-term variation of the asymptotic series
is justif\/ied (see the remarks below),
\begin{equation}\label{cix}
\var{a_k}[A]=(n-k)a_k[\w,A].
\end{equation}

From the zeta function point of view, the variational calculation is
\begin{gather*}
\G(s)\var{\z(s)}=\G(s)(\Tr\,A^{-s})^\bul =\G(s)(\Tr\,D^{-hs})^\bul =\G(s)\Tr\,\var{D}(-hs)D^{-hs-1} \\
\phantom{\G(s)\var{\z(s)}}{}=(b-a)hs\G(s)\Tr\,\w D^{-hs} =2\ell s\G(s)\z(\w,s).
\end{gather*}
But $a_k$ is the residue of the $\G\z$ function at $(n-k)/2\ell$,
so we recover~\nnn{cix}.

In addition, assuming for the moment that 0 is not an eigenvalue of $A$,
\[
\var{\z'_A(0)}=2\ell\left.\dfrac{d}{ds}\right|_{s=0}(s\z_A(\w,s))
=2\ell\z_A(\w,0).
\]
But $\z_A(\w,0)=a_n[\w,A]$, since ${\rm Res}_{s=0}\G(s)=1$.
If 0 is an eigenvalue for $A$, we carry out the last calculations
after replacing $\exp(-tA)$ and
$A^{-s}$ by the restrictions of $\exp(-t|A|)$ and $|A|^{-s}$ to
the orthogonal complement in $L^2$ of $\cN(A)$, the null space of $A$.
(Here a ``restriction'' is extended back to the whole space as
0 on $\cN(A)$.)
On the level of kernel functions, we are then dealing with
\[
\displaystyle\sum_{\l_j\ne 0}\exp(-|\l_j|t)\f_j\otimes\f_j^*
\qquad{\rm and}\qquad
\sum_{\l_j\ne 0}|\l_j|^{-s}\f_j\otimes\f_j^*.
\]
The zeta function arises from such a modif\/ied
kernel, while the $Z$-function (as ref\/lected in the~$U_k[A]$) arises
from the ``original'' kernel function.  Chasing through the consequences
of this mismatch gives rise to the extra term in~\nnn{vardet}.

In both the $Z$-function and zeta function arguments for the
variational formula \nnn{cit}, there is an interchange of limit
operations.  In the case of the $Z$ function, we are claiming that a
certain asymptotic series may be dif\/ferentiated term-by-term with
respect to an outside parameter (the $\e$ in the conformal factor
$e^{\e\w}$).  That is, we are asserting that the asymptotics of the
variation are the variation of the asymptotics.  In the zeta function
argument, an equivalent strong assertion is made: that the analytic
continuation of the variation (again, with respect to the outside
parameter $\e$) is the variation of the analytic continuation.  This
is not just a matter of the uniqueness of the analytic continuation.
In general, families of holomorphic functions with extremely nice
dependence on an outside parameter may exhibit very ugly dependence in
their analytic continuations.  In either case, the issue comes down to
suitable uniformity in the outside parameter of the estimates of error
in the asymptotic expansion -- in the zeta function
argument, this enters when
the heat operator and the Mellin transform are used to
accomplish the analytic continuation.

It's also worth noting that the issue of term-by-term dif\/ferentiation
of the asymptotic series with respect to an outside parameter (here,
$\e$) is more delicate than term-by-term dif\/ferentiation with
respect to the asymptotic parameter $t$.  Indeed, since one can
trivially {\em integrate} term by term, one can also dif\/ferentiate: if
\begin{gather*}
\displaystyle f(t)\sim\sum_ka_kt^k, \qquad
f'(t)\sim\sum_kb_kt^k\qquad {\rm with}\quad b_0=0,
\end{gather*}
then
\[
f(t)+{\rm const}\sim\sum_k\dfrac{b_k}{k+1}t^{k+1},
\]
so that
\[
a_{k+1}=\dfrac{b_k}{k+1},\qquad k\ne -1.
\]
This shows that the series for $f'(t)$ is the term-by-term derivative
of that for $f(t)$, since the derivative of the term $a_0t^0$ is 0.

Note that the determinant quotient problem for conformal metrics is
more or less trivial in odd dimensions (at least for manifolds without
boundary), since $U_n[A]=0$.  (As noted earlier, the $a_{{\rm odd}}$
vanish; the stronger statement that the $U_{{\rm odd}}$ vanish is also
true.)  This means that the only variation in the determinant comes
from the $\f_j$ terms in~\nnn{vardet}.

Here's a sample calculation to illustrate the dimensional regularization
method of computing the determinant quotient.  Our ``input'' operator
$A$ will be the square of the Dirac operator $\dc=\g^a\nd_a$
on the bundle of spinors.  (Recall the discussion just before the Theorem
on Variations.)  Suppose we would like the determinant quotient in
dimension 2.  Since
\[
a_2[\dc^2]=-\dfrac1{12}\,2^{[n/2]}(4\pi)^{-n/2}\int K,
\]
the determinant quotient will be
\begin{equation}\label{dcsample}
\left[\dfrac1{n-2}\dint((K\,dv)_\w-(K\,dv)_0)\right]_{n=2}.
\end{equation}
If this is to work, there must somehow be a ``hidden'' $n-2$ factor
in the integral.
To bring this out, recall the {\em Yamabe
equation}
\begin{equation}\label{ybeeq}
Yu:=\D u+\dfrac{n-2}{4(n-1)}Ku=
\dfrac{n-2}{4(n-1)}\overline{K}u^{\frac{n+2}{n-2}}
\end{equation}
for the scalar curvature after a conformal change
\[
\bar g=\W^2g,\qquad 0<\W\in\ci,\qquad u=\W^{\frac{n-2}2}.
\]
(Equation \nnn{ybeeq}
is actually the conformal covariance relation for the conformal
Laplacian, applied to the function~$1$.)
Rewrite \nnn{ybeeq} as
\begin{gather*}
\D\big(e^{\frac{n-2}2\w}\big)
+\dfrac{n-2}{4(n-1)}Ke^{\frac{n-2}2\w}=
\dfrac{n-2}{4(n-1)}\overline{K}e^{\frac{n+2}2\w}, \\
\D_0\big(e^{\frac{n-2}2\w}-1\big)
+\dfrac{n-2}{4(n-1)}K_0e^{\frac{n-2}2\w}
=\dfrac{n-2}{4(n-1)}K_\w e^{\frac{n+2}2\w}.
\end{gather*}
On the second line, we have introduced a gratuitous $-1$ into the
argument of $\D_0$; this is harmless because $\D_01=0$.
It is useful because it allows us to
divide by $(n-2)/2$:
\begin{equation}\label{glean}
\D_0f\left(\frac{n-2}2\w\right)
+\dfrac{1}{2(n-1)}K_0e^{\frac{n-2}2\w}
=\dfrac{1}{2(n-1)}K_\w e^{\frac{n+2}2\w},
\end{equation}
where
\[
f(x):=\dfrac{e^x-1}{x}.
\]
Evaluation at $n=2$ now yields the
{\em Gauss curvature prescription equation}
\[
\D_0\w+\dfrac12K_0=\dfrac12K_\w e^{2\w}\qquad (n=2).
\]
In particular,
\begin{equation}\label{intKdc}
\int(K\,dv)_0=\int(K\,dv)_\w\qquad (n=2),
\end{equation}
since the dif\/ference of the two integrands, $2(\D_0\w)dv_0$, integrates
to zero by the Divergence Theorem.

\nnn{intKdc} just expresses the invariance of the conformal index
in this setting.  To get a handle on~\nnn{dcsample} we need to
go to higher order in the above calculation, in the dimension parameter
\[
\b=\dfrac{n-2}2.
\]
For this, f\/irst
multiply \nnn{glean} by $dv_0=e^{-n\w}dv_\w$:
\[
\left\{\D_0\left(\dfrac{e^{\b\w}-1}{\b}\right)
+\dfrac1{2(n-1)}K_0e^{\b\w}\right\}dv_0
=\dfrac1{2(n-1)}K_\w e^{-\b\w}dv_\w\,,
\]
Now multiply by $\exp(\b\w)$:
\[
\left\{e^{\b\w}\D_0\left(\dfrac{e^{\b\w}-1}{\b}\right)
+\dfrac1{2(n-1)}K_0e^{2\b\w}\right\}dv_0
=\dfrac1{2(n-1)}(K\,dv)_\w.
\]
If we now set
\[
J:=\dfrac{K}{2(n-1)},
\]
the above can be rewritten as
\begin{gather*}
(J\,dv)_\w-(J\,dv)_0=
e^{\b\w}\D_0\left(\frac{e^{\b\w}-1}{\b}\right)
+(e^{2\b\w}-1)(J\,dv)_0 \\
\phantom{(J\,dv)_\w-(J\,dv)_0}{}=\D_0\left(\frac{e^{\b\w}-1}{\b}\right)
+(e^{\b\w}-1)\D_0\left(\frac{e^{\b\w}-1}{\b}\right)
+(e^{2\b\w}-1)(J\,dv)_0.
\end{gather*}

Now integrate and divide by $n-2$:
\[
\dfrac1{2\b}\dint\left((J\,dv)_\w-(J\,dv)_0\right)
=\dint\left\{\dfrac{e^{\b\w}-1}{2\b}\D_0
\left(\dfrac{e^{\b\w}-1}{\b}\right)
+\dfrac{e^{2\b\w}-1}{2\b}J_0\right\}dv_0\,.
\]
Evaluating at $n=2$ (i.e., at $\b=0$), this gives
\[
\left[\dfrac1{n-2}\dint\left((J\,dv)_\w-(J\,dv)_0\right)\right]_{n=2}
=\dint\left\{\dfrac12\w\D_0\w+\w J_0\right\}dv_0.
\]
As a result, if the Dirac operator does not take the eigenvalue 0
in the conformal class in which we are working,
\[
\z'_{\minidcs}(0)
-\z'_{\minidcso}(0)=
-\dfrac1{12\pi}\dint\left\{\dfrac12\w\D_0\w+\w J_0\right\}dv_0.
\]

Note that the objects that appear in the f\/inal result,
$\D$ and $J$, really appear just
because of the form of the integrated invariant $a_2$.
In fact, since the only local scalar invariant
that can contribute to $a_2$ is $K$ (or $J$), the above formula
is independent of $A$, except for the overall constant factor.
For example, for the Laplacian on functions in dimension~2,
\[
\z'_{\D_\w}(0)-\z'_{\D_0}(0)=\dfrac1{12\pi}
\dint\left\{\dfrac12\w\D_0\w+\w J_0\right\}dv_0.
\]
We compute this by extending $\D$ for $n=2$ to the conformal Laplacian
$\D+(n-2)J/2$ in higher dimensions.

\section{Extremals of the functional determinant}

When we begin to think about maximizing or minimizing
the determinant, the f\/irst thing we notice is
that it is not invariant
under uniform scaling of the metric:
\[
\tilde g=e^{2\a}g,\quad \a\in\bbR\quad \Rightarrow\quad
\tilde A=e^{-2\ell\a}A,
\]
so
\begin{gather*}
\z_{\tilde A}(s)=e^{2\ell\a s}\z_A(s), \qquad
\z'_{\tilde A}(0)=\z'_A(0)+2\ell\a\z_A(0).
\end{gather*}
Thus it's necessary to somehow {\em penalize} the determinant for
the size of the metric.  One way to do this is with a volume penalty:
\[
\cD(A,g):=({\rm vol}\,g)^{2\ell\z_A(0)/n}\det A_g
\]
is invariant under uniform dilation.
(Recall that $\z_A(0)$ is a conformal invariant.)
The formula for the conformal variation of the determinant
says that the conformal variation of the $\cD$-func\-tional~is
\[
\var{-\log(\cD_A)}=2\ell{\rm vol}(g)\big((\w U_n)\tilde{\ }-
\tilde\w\tilde U_n\big),
\]
where the tilde denotes the average value.  Thus the $\cD$-quotient
is
\[
\log\dfrac{\cD(A_\w)}{\cD(A_0)}=
\dfrac{2\ell\z_A(0)}{n}\log
\underbrace{\frac{\int e^{n\w}dv_0}{\int dv_0}}_{(e^{n\w})\tilde{\ }}
+\log\dfrac{\det(A_\w)}{\det(A_0)}.
\]

Now (to be somewhat vague) the coef\/f\/icients in $\z_A(0)$
{\em have something to do with} the coef\/f\/i\-cients in the determinant
term -- $\z_A(0)$ is $a_n[A]$, and the determinant term somehow
{\em came from}~$a_n[A]$.
As a result, though dif\/ferent $A$ produce dif\/ferent $\cD$-functionals,
combinations of a~limited list of functionals keep appearing.
For example in dimension~2, in the standard conformal class
on the sphere ${\bf S}^2$, one is always looking at
\[
-\log\int_{\bbS^2}e^{2(\w-\bar\w)}d\x+\int_{\bbS^2}\w(\D_0\w)d\x,
\]
where $d\x$ is the standard measure, normalized so that $\int_{{\bf S}^2}d\x=1$.
This is {\em exactly} the quantity that appears in the celebrated
{\em Moser--Trudinger
inequality}.

Recall the {\em Sobolev embeddings}
\[
L^2_r\hookrightarrow L^{2n/(n-2r)}.
\]
The inequalities asserting the boundedness of the inclusion operators
have {\em sharp forms} (more on this later; see \nnn{beclp} below) on ${\bf S}^2$.
What happens as $r\to n/2$?  One answer is the {\em Sobolev lemma}
\[
L^2_{\frac{n}2+\e}\hookrightarrow C^0.
\]
A deeper answer is the embedding
\[
L^2_{n/2}\hookrightarrow e^L.
\]
The inequality describing this inclusion
compares an $L^2_{n/2}$ norm with a
functional that governs admission to the {\em Orlicz class} $e^L$.
When $n=2$, this is the Moser--Trudinger inequality;
for other~$n$, it is the exponential class inequality of Beckner
\cite{bec} and Carlen--Loss \cite{cl}.  The sharp form of the
Moser--Trudinger inequality is:

\begin{theorem}[Moser--Trudinger inequality]
\[
\log\int_{\bbS^2}e^{2(\w-\bar\w)}d\x\le\int_{\bbS^2}\w(\D_0\w)d\x,
\]
with equality iff $e^\w$ is a constant multiple of a conformal
factor $\W_h$ on $\bbS^2$.
\end{theorem}

To explain the term {\em conformal factor}
as it's being used here, recall that there is a
{\em M\"obius group} of
conformal transformations of the Riemann sphere.  In terms of
Riemannian geometry, such a transformation (dif\/feomorphism) has
\begin{equation}\label{conftr}
(h^{-1})^*g_0=\W_h^2g_0.
\end{equation}
The inequality is {\em sharp} because equality is attained for some
$\w$, and the case of equality is completely analyzed.

In 1982, Onofri \cite{ono} used the Moser--Trudinger inequality to
prove that on ${\bf S}^2$, the volume-penalized determinant of $\D$ is
maximized exactly at the standard metric and its dif\/feomorphic images.
Part of this is special to 2 dimensions~-- each metric on ${\bf S}^2$ is
dif\/feomorphic to one in the standard conformal class.  Since
dif\/feomorphisms preserve spectral invariants, it's therefore
enough to study the
standard conformal class.  This is not the case in higher dimensions,
even on the sphere.

In fact, the space which may be loosely described by
\[
\dfrac{{\rm metrics}}{{\rm Diffeo}(M)
\ltimes
{\rm Conf}(M)}
\]
is just a point when $M={\bf S}^2$.  (The semidirect product $\ltimes$ in
the denominator ref\/lects the fact that dif\/feomorphisms act on
conformal factors.)
This space is generally a f\/inite-dimensional
object when $n=2$ (in fact, {\em Teichm\"uller space}), but
it's a big, wild object when $n>2$.

One of the many questions that one can ask at this point is
how all this plays out in higher dimensions~-- on ${\bf S}^4$, ${\bf S}^6$, and
so on.  Can we still solve a max/min problem by using sharp
inequalties?  What will
the $\cD$-quotient functional and the sharp inequalities
look like?  In the following paragraphs, we look at the objects and
methods that will generalize the two-dimensional
versions appearing in the paragraph headings.

{\bf The Laplacian.}
The Laplacian $\D$ on scalar functions
turned up in the formula for the determinant quotient
in dimension 2,
even if it was some ``exotic'' $A$ (like Dirac-squared)
whose determinant we were measuring.
In 1992, Graham, Jenne, Mason, and Sparling \cite{gjms} proved the existence
of a series of conformally covariant operators $P_m$,
of order $m$ and of the form $\D^{m/2}+{\rm(lower\ order)}$;
these operators exist for
\begin{gather*}
{\rm all\ even}\ m\ {\rm when}\ n\ {\rm odd}, \\
{\rm even}\ m\le n\ {\rm when}\ n\ {\rm even}.
\end{gather*}
The {\em Paneitz operator} $P_4$ was previously found in
\cite{pan83,rieg84,es85}.

In \cite{wunsch86}, W\"unsch wrote a kind
of formula for $P_6$ -- really more of an algorithm for computing it,
but enough to prove its existence.

The role of the 2-dimensional Laplacian in the determinant formulas is
played by $P_n$ (the operator of dimension order), for even $n$.
The role of the Gauss curvature is played by an invariant
$Q_n$ (discovered in \cite{bo91} in dimension 4, and in
\cite{tbkorea} in higher even dimensions) for which
\begin{equation}\label{qm}
P_m=P_m^0+\dfrac{n-m}2Q_m,\qquad P_m^01=0.
\end{equation}

{\bf Dimensional regularization.}
The algebra of local invariants for general metrics can be simplif\/ied
somewhat in the category of {\em conformally flat} metrics.  These
are characterized by the existence of local coordinates near each point
in which
the metric is $g_{ab}=\W^2\d_{ab}$ for some smooth positive
function $\W$.  There is an alternative characterization, in terms
of invariants:
if $n\ge 4$, a metric is conformally f\/lat if and only if its {\em Weyl
tensor} vanishes.
The Weyl tensor, in turn, may be characterized as the totally trace-free
part of the Riemann curvature tensor.  One has the following formulas.
Let
\[
J=\frac{K}{2(n-1)}
\]
as above, and let
\[
V=\frac{r-Jg}{n-2},
\]
where $r$ is the Ricci tensor.
Then the Weyl tensor $C$ is
\[
C_{abcd}=R_{abcd}+V_{bc}g_{ad}-V_{bd}g_{ac}+V_{ad}g_{bc}-V_{ac}g_{bd}.
\]
As a consequence of the Bianchi identity, we have
\begin{equation}\label{xtra}
\N^aC_{abcd}=(n-3)(\N_cV_{bd}-\N_dV_{bc}).
\end{equation}
When $n=3$, the vanishing of $\N_cV_{bd}-\N_dV_{bc}$ characterizes
conformal f\/latness (a fact which in itself is suggestive of some kind
of dimensional regularization).  In any case, for $n\ge 3$, in the
conformally f\/lat case, all local invariants are built from $g$, $V$ and
$\N$ ($J$ being an abbreviation for $V^a{}_a$),
and one has the ``extra identity'' that $\N V$ is totally symmetric.
The ideal of identities (i.e., the ideal of left sides of identities with
right side~0) is generated by~\nnn{xtra}, together with the Ricci
identities, which give formulas for $[\N_a,\N_b]$ applied to tensors
of all ranks.

We shall concentrate on dimensional regularization in the conformally
f\/lat category, though the restriction is not essential~-- there
is a perfectly good algebraic setting in the category of arbitrary
metrics.  However, conformally f\/lat metrics provide
the cleanest route to the
very clean extremal results on the
spheres
${\bf S}^{{\rm even}}$.

Recall now our earlier formula for the determinant quotient:
\[
\dfrac1{2\ell}\left(\z'_{A_\w}(0)-\z'_{A_0}(0)\right)
=\left.\dfrac{a_{n}[A_\w]-a_{n}[A_0]}{n'-n}\right|_{n'=n},
\]
where $n$ is the f\/ixed dimension (in which we want a formula),
and $n'$ is the ``running'' dimension.  What one gets out of this
for even $n>2$, given suitable dependence on $n'$, is a formula
for the penalized determinant quotient
functional of the form
\begin{gather}\label{dfcnlform}
C\left(\frac{n}{2(Q_n)_0}\int_{\bbS^n}(P_n\w,\w)d\x
-\log\int_{\bbS^n}e^{n(\w-\bar\w)}d\x\right)
+\int_{\bbS^n}\left((B\,d\x)_\w-(B\,d\x)_0\right),
\end{gather}
where $B=B[A]$ is some local scalar invariant, and $C$ is a constant.

The invariant $B$ is a linear combination
of invariants from a precise (and remarkably short, for small $n$)
list.
For $n=2$, this list is empty; for $n=4$, there is one invariant that
can contribute;
for $n=6$, there are 3, and for $n=8$, there are 8.

\begin{remark}
Besides being a convenient tool for calculation, dimensional regularization
produces a theoretical result on the nature of the functional, which
may be paraphrased as follows.
If one follows the other, f\/ixed-dimension method of integrating along
a curve in the conformal class, then after subtracting out the
correct multiple of the functional
\[
\cF_0(\w):=
\dfrac{n}{2(Q_n)_0}\int_{\bbS^n}(P_n\w,\w)d\x
-\log\dint_{\bbS^n}e^{n(\w-\bar\w)}d\x,
\]
the dependence on $\w$ of the remainder is that of a local invariant:
\[
\cF(\w):=\dint_{\bbS^n}\left((B\,d\x)_\w-(B\,d\x)_0\right).
\]
\end{remark}

We saw that in dimension 2, the penalized determinant quotient
functional was an exactly the one estimated by
the Moser--Trudinger inequality.  What is the analogue for dimension
4, and higher dimensions?  Remarkably, in the early 1990's,
several lines of thought necessary for an understanding of this
were maturing -- higher dimensional Polyakov formulas
\cite{bo91}, sharp inequalities \cite{bec,cl}, and,
though only an understanding of the Paneitz operator is needed
in 4~dimensions, the GJMS operators \cite{gjms}.

The sharp exponential class inequality in $n$ dimensions is:

\begin{theorem}[Exponential class inequality]
The quantity
\begin{equation}\label{expocls}
\dfrac{n}{2\G(n)}\int_{\bbS^n}((P_n)_0\w,\w)d\x
-\log\int_{\bbS^n}e^{n(\w-\bar\w)}d\x
\end{equation}
is $\ge 0$, with $=$ iff $e^\w$ is a conformal factor.
\end{theorem}

The term {\em conformal factor} is again meant in the sense of
\nnn{conftr}.  The analogue of the M\"obius group, the conformal
group of ${\bf S}^n$, is isomorphic to
SO$(n+1,1)$; the action of this group will be needed in more detail
later.

The operator $(P_n)_0$ is the GJMS operator
evaluated at the round metric; the precise form of this operator
was given in \cite[Remark 2.23]{tbjfa87}.  Namely,
let
\[
\cB=\sqrt{\D+\left(\dfrac{n-1}2\right)^2}\,.
\]
Then
\begin{equation}\label{old87}
(P_m)_0=\prod_{a=1}^{m/2}\underbrace{\left\{
\left(\cB+a-\frac12\right)\left(\cB-a+\frac12\right)
\right\}}_{\cB^2-(a-\frac12)^2}.
\end{equation}

Recall that the GJMS operator implicitly def\/ined a local invariant
via \nnn{qm}.  Let
\[
\dfrac{n-m}2Q_m=P_m1,
\qquad \mbox{so that} \qquad
P_m=P_m^0+\dfrac{n-m}2Q_m,
\]
where $P_m^0=1$.  In general, $P_m$ and $Q_m$ are
not unique -- for example,
one may add any multiple of the norm-squared of the Weyl tensor,
$|C|^2=C^{abcd}C_{abcd}$, to the Paneitz operator $P_4$ and still
have a conformal covariant.  What one really needs, for present
purposes, is a choice of $P_m$ that~is
\begin{itemize}
\itemsep=0pt
\item has principal part $\D^{m/2}$;
\item has conformal biweight $((n-m)/2,(n+m)/2)$;
\item has coef\/f\/icients which are rational in the dimension $n$;
\item is formally self-adjoint; and
\item annihilates the function 1 in dimension $n$.
\end{itemize}

Of course, such a $P_m$ is neither an operator on a particular
manifold, nor an operator scheme for Riemannian manifolds in
a certain dimension, but rather a {\em formula}, involving the
dimension and some basis of the possible terms (built from
$g$, $\N$, and $R$) which is stable for large dimensions.

There is certainly such a $P_m$
in the category of conformally f\/lat metrics, and there is even
an algorithm for writing it.  Since such an operator must
be $\D^{m/2}$ for the f\/lat metric $g_{ab}=\d_{ab}$, we just
need to write the operator invariantly for metrics
$(g_\w)_{ab}=e^{2\w}\d_{ab}$.  By the conformal covariance relation,
\[
(P_m)_\w=e^{-(n+m)\w/2}(P_m)_0e^{(n-m)\w/2}.
\]
(As usual, we blur the distinction between a function $f$ and
the operator $\m_f$ of multiplication by~$f$~-- the f\/inal exponential
in the foregoing formula should really be a multiplication operator.)
In order to write this in terms of covariant derivatives in the
metric $g_\w$, we write $(P_m)_0$ as a~composition of conformal
covariants evaluated at the metric $g_0$, and use the conformal
covariance relations for these ``smaller'' operators.  For
example,
\[
(P_m)_0=\D_0^{m/2}=Y_0^{m/2}=\big(e^{(n+2)\w/2}Y_\w
e^{-(n-2)\w/2}\big)^{m/2},
\]
where $Y$ is the conformal Laplacian.
As a result,
\[
(P_m)_\w=e^{-(n+m)\w/2}
\big(e^{(n+2)\w/2}Y_\w
e^{-(n-2)\w/2}\big)^{m/2}
e^{(n-m)\w/2}.
\]

Another way to compute would be to use the fact that $\D=\d d$,
and write
\[
(P_m)_0=(\d_0d)^{m/2}=\big(e^{n\w}\d_\w e^{-(n-2)\w}d\big)^{m/2},
\]
using the metric invariance of $d$, and the conformal biweight
$(n-2,n)$ of $\d$.  The result is
\[
(P_m)_\w=e^{-(n+m)\w/2}
\big(e^{n\w}\d_\w e^{-(n-2)\w}d\big)^{m/2}
e^{(n-m)\w/2}.
\]

At any rate, we now commute all covariant
derivatives to the right, past the expressions in~$\w$.  This produces
terms in covariant derivatives of $\w$.  However, all second and
higher derivatives of $\w$ may be eliminated, in favor of local
invariants, by iterated use of the Ricci tensor change formula
\[
\w_{ab}=-V_{ab}-\w_a\w_b+\frac12\w_c\w^cg_{ab}
\]
and its iterated covariant derivatives.  (Here
$\w_{a\ldots b}$ abbreviates $\N_b\cdots\N_a\w$.)  The total degree of each
term in $e^\w$ is 0, so we only have to worry about leftover occurrences
of the f\/irst derivative,~$\N\w$.
There are none such, however, as a consequence of our
having picked the correct principal part and biweight for the conformal
covariant $P_m$.
The result is an expression for $P_m$ in the conformally f\/lat
category.

In the case of high-order bundle-valued conformal covariants,
a procedure like that just above is also possible, given the correct
principal part and biweight.  (The correct principal part
may be deduced from spectral considerations \cite{tbsw,ab2}.)
In general, the principal part
is written as a linear combination of compositions of ``smaller''
conformal covariants.

The subject of {\em promoting} conformal covariants in the conformally f\/lat
category to conformal covariants
for general metrics still features some open questions.  It is known
\cite{graham} that the operator on scalar functions with principal
part $\D^3$ does not promote in dimension 4; note that this
is just beyond the restriction $m\ge n$ on the GJMS operators in
even dimensions.  In terms of Bernstein--Bernstein--Gelfand resolutions
(see, e.g., \cite{eastslov}), which classify (dif\/ferential)
conformal covariants
in the conformally f\/lat category, one knows that operators which
are not {\em longest arrows} do in fact promote.
Some longest arrows (like the critical order GJMS operator~$P_n$)
promote as well.  A plausible conjecture is that the $P_n$ are the
{\em only} longest arrows that promote.

If we stay in the conformally f\/lat category it is easily shown,
by examining the above construction (which produces
unique $P_m$), that the
above conditions on $P_m$ hold.
In addition \cite{sharp},
with the subscript 0 still denoting evaluation at
the round sphere metric,
\[
(Q_m)_0=\dfrac{\G\left(\frac{n+m}2\right)}
{\G\left(\frac{n-m+2}2\right)}, \qquad
{\rm and\ in\ particular},\qquad (Q_n)_0=\G(n).
\]
This means that the f\/irst functional in
\nnn{dfcnlform}, which we have called $\cF_0(\w)$, is exactly the
quantity estimated (i.e.\ asserted to be positive) by the exponential
class inequality:
\[
\cF_0(\w):=
\dfrac{n}{2\G(n)}\int_{\bbS^n}(P_n\w,\w)d\x
-\log\dint_{\bbS^n}e^{n(\w-\bar\w)}d\x.
\]

The issue is now whether we can get an estimate of the other term,
\[
\cF(\w):=\dint_{\bbS^n}\left((B\,d\x)_\w-(B\,d\x)_0\right),
\]
which is
{\em compatible} with the exponential class estimate.
This means that it should assert that
\[
C\cF(\w)\ge 0,
\]
where $C$ is the constant in \nnn{dfcnlform} -- the inequality should
``go in the same direction as'' the exponential class one.
That is, it should assert this provided $C\ne 0$.  If $C$ vanishes,
then $\z_A(0)$ vanishes; one way to see this is that the log term
comes from the volume penalty; the case when no penalty is needed is
exactly $\z_A(0)=0$.  It is not clear whether this can happen in
reality; it might be reasonable to conjecture that $\z_A(0)\ne 0$
for all $A$ satisfying our assumptions.

\begin{remark} It would also be nice if
our compatible inequalities were compatible in another way~--
their cases of equality should contain the case of equality for
the exponential class inequality.  This will assure that the
coupled functional vanishes exactly at conformal factors.
But this additional property is provided automatically by the
invariant nature of $\cF(\w)$ (Remark~3)~-- the
functional is dif\/feomorphism-invariant, and thus in particular
is invariant under conformal dif\/feomorphisms.  This guarantees
that when $\w$ is a conformal factor, $\cF(\w)$ will vanish.
\end{remark}

To be more specif\/ic on how $P_n$ entered our formula for
the $\cF_0(\w)$ part of the determinant quotient, recall f\/irst
how the Gauss curvature
prescription equation arose from the Yamabe equation.
More generally, we would like a $Q_n$ prescription problem to arise
in the same way from the $Q_m$ prescription problems, by dimensional
regularization:
\begin{gather*}
P_m=P_m^0+\dfrac{n-m}2Q_m, \qquad P_m^01=0, \qquad
e^{\frac{n+m}2\w}(P_m)_\w f=(P_m)_0\big(e^{\frac{n-m}2\w}f\big), \\
\dfrac{n-m}2(Q_m)_\w e^{\frac{n+m}2\w}=\left(P_m^0+\dfrac{n-m}2Q_m\right)_0
e^{\frac{n-m}2\w} \\
\phantom{\dfrac{n-m}2(Q_m)_\w e^{\frac{n+m}2\w}}{}
=(P_m^0)_0\big(e^{\frac{n-m}2\w}-1\big)+\dfrac{n-m}2\left(Q_m\right)_0
e^{\frac{n-m}2\w}.
\end{gather*}
Dividing by $(n-m)/2$ and evaluating at $n=m$, we have
\[
(Q_n)_\w e^{n\w}=(P_n^0)_0\w+(Q_m)_0.
\]
In particular, since $P_n^0$ has the form $\d Sd$ (by some invariant theory),
\[
\int(Q_ndv)_\w=\int(Q_ndv)_0.
\]

Going up an additional order in $\b=(n-m)/2$, as in the $m=2$ case, we get
\[
\dint\left.\dfrac{(Q_ndv)_\w-(Q_ndv)_0}{n'-n}\right|_{n'=n}=
\frac12\int\w(P_n)_0\w\,dv_0
+\dint\w(Q_n)_0dv_0\,.
\]
In this way, $P_n$ has naturally entered the determinant calculation.

Here's the situation in dimension $n=4$.  $P_4$ is the {\em Paneitz
operator}:
\[
P_4=\D^2+\d Td+\dfrac{n-4}2Q_4,
\]
where
\[
T=(n-2)J-4V\cdot, \qquad
(V\cdot\f)_b=V^a{}_b\f_a, \qquad
Q_4=\dfrac{n}2J^2-2|V|^2+\D J.
\]

Since $\int U_4[A]$ is a conformal invariant, $U_4[A]$, which {\em a priori}
is a linear combination of $K^2$, $|r|^2$, $|R|^2$, and $\D K$
(or, equivalently, $J^2$, $|V|^2$, $|C|^2$, and $\D J$), must
be something more special, namely a linear combination
of $Q_4\,$, $\D J$, and $|C|^2$ ($C={\rm Weyl}$).
Since these quantities have conformally invariant integrals,
the leftover functional
\[
\int_{\bbS^4}\left((B\,d\x)_\w-(B\,d\x)_0\right)d\x
\]
can be expressed
as a multiple of
\[
\int_{\bbS^4}\left((J^2d\x)_\w-(J^2d\x)_0\right)d\x.
\]
($J^2$ could actually be replaced in this formula by anything
linearly independent of $Q_4\,$, $\D J$, and~$|C|^2$.)

In higher dimensions, because of the Conformal Index Theorem,
there is always some
leeway in the choice in $B$.
By \cite{bgp95}, if $L$ is a local invariant with $\int L$ conformally
invariant, then in the conformally f\/lat category,
\[
L={\rm const}\cdot{\rm Pf\/f}+{\rm div},
\]
where ``div'' is an exact divergence, and ``Pf\/f''
is the Pfaf\/f\/ian (Euler characteristic density).
Since $\int Q_n$ and $\int U_n$ are conformally invariant,
\[
U_n[A]={\rm const}\cdot Q_n+{\rm div},
\]
where the constant is uniquely determined by $A$.

In dimension 4, the only possible universal divergence is $\D J$.
Thus $J^2$ appears, basically, because it is a conformal primitive
for $\D J$:
\[
\left(\int J^2\right)^{\var{\ }}=2\int J\D\w=2\int(\D J)\w.
\]

As an overall conclusion in dimension 4, all of
our penalized determinant quotient functionals end up being linear
combinations of the two functionals
\begin{gather*}
\cF_0(\w)=\frac13\dint_{\bbS^4}\w\D_0(\D_0+2)\w\,d\x-\log\dint_{\bbS^4}
e^{4(\w-\bar\w)}d\x, \\
\cF_1(\w)=\dint_{\bbS^4}\{(J^2d\x)_\w-(J^2dx)_0\}.
\end{gather*}

By the above, we know all about the functional $\cF_0(\w)$.
The other functional,
$\cF_1(\w)$, is conveniently analyzed by
the solution of the Yamabe
problem \cite{ybe,trud,aubin,schoen}.
The {\em Yamabe quotient} of the sphere $\bbS^4$ is
\[
\cY(f)=\dfrac{\dint_{\bbS^4}f(\D_0+2)f\,d\x}
{\left(\dint_{\bbS^4}f^4d\x\right)^{1/2}}.
\]
The minimum value of $\cY(f)$ is attained at $f=1$, and so is $2$.
Thus for nowhere-vanishing $f$,
\begin{equation}\label{hold}
2\|f\|_4^2\le\int_{\bbS^4}\underbrace{f(\D_0+2)f}_{f^2\cdot\left(
\frac{\D_0f}f+2\right)}
\,d\x
\le\underbrace{\|f^2\|_2}_{\|f\|_4^2}\left\|\frac{\D_0f}f+2\right\|_2,
\end{equation}
by the Schwarz inequality, so that
\[
2\le\left\|\frac{\D_0f}f+2\right\|_2.
\]
We would like to apply this to the nowhere-vanishing function
\[
f=e^{\frac{n-2}2\w}=e^\w.
\]
This yields
\[
\displaystyle 4\le\dint_{\bbS^4}\Bigg(e^{-\w}(
\underbrace{\D_0+2}_{Y_0})e^\w
\Bigg)^2d\x
=\sdint_{\bbS^4}J_\w^2dv_\w.
\]
This result is well worth stating in its own right:

\begin{theorem}
$
0\le\dint_{\bbS^4}\left((J^2d\x)_\w-(J^2dx)_0\right),
$
with $=$ iff $e^\w$ is a conformal factor.
\end{theorem}

The last part, about the case of equality, comes from the case of equality
for the Yamabe problem on $\bbS^4$.
The inequality really expresses the instance
$L^2_1\hookrightarrow L^4$ of the Sobolev Embedding Theorem in
a way adapted to spectral invariant theory.
As it has turned out, the functional $\cF_1(\w)$ has
{\em exactly} the same case of equality as the exponential class
functional $\cF_0(\w)$, not just a containing case of equality as
described in Remark~4.

Just what are the coupling constants between the functionals
$\cF_0(\w)$ and $\cF_1(\w)$ for some prominent admissible operators
$A$?  Up to overall constant factors, here are the functionals when
$A$ is the conformal Laplacian, Dirac-squared, and Paneitz operators:
\begin{alignat*}{3}
&Y:&&-3\cF_0(\w)-2\cF_1(\w), &\\
&\dc^2:\ \ &&33\cF_0(\w)+7\cF_1(\w),& \\
&P_4:&& 21\cF_0(\w)-16\cF_1(\w).&
\end{alignat*}
Note that by dimensional regularization, these numbers really fall
out of formulas for $a_4[A]$, provided we keep enough information
to read the dimension dependence.  (See~\cite{tbcmp96}
for the $P_4$ calculation.)

Thus the extremal problem for $\cD(Y)$ and $\cD(\dc^2)$ is solved,
in the standard conformal class on $\bbS^4$.  This calculation was
originally done in \cite{bcy}, based on the determinant formula
in \cite{bo91}.

For the penalized $\det\,P_4$, functional, it's still unclear
whether the round metrics are absolute extremals.
The question can be phrased as comparing the sizes
of the gaps (the dif\/ferences between the left and right sides)
in two sharp inequalities.
The {\em leading terms} (those with the highest number of
derivatives of $\w$)
for $\cF_0+a\cF_1$ produce a quadratic form that is
\begin{alignat*}{3}
&\mbox{{\rm positive\ def\/inite\ for}}\ \ &&\  a\le -8/15,& \\
&\mbox{{\rm indef\/inite\ for}}\ \ && \ -8/15 <a<-1/3, &\\
&\mbox{{\rm negative\ def\/inite\ for}}\ \ && \ a\ge-1/3.&
\end{alignat*}
The number for $P_4$ is $-16/21<-8/15$; this is suggestive
of a minimum for $\cD((P_4)_\w)$ at the round metrics,
but
the question is still open.

A look at dimension 6 indicates that something interesting
is going on.
The $\cD$-functional in the standard conformal class on $\bbS^6$
is a linear combination of four functionals,
\begin{gather}
\cF_0(\w)=\tfrac1{40}\dint_{\bbS^6}\w((P_6)_0\w)d\x-\log\dint_{\bbS^6}
e^{6(\w-\bar\w)}d\x, \qquad
\cF_1(\w)=\dint_{\bbS^6}(|dJ|^2d\x)_\w, \nonumber\\
\cF_2(\w)=\dint_{\bbS^6}((|dJ|^2+2J^3)d\x)_\w-54, \nonumber\\
\cF_3(\w)=\dint_{\bbS^6}\left(\left(|dJ|^2+\tfrac{28}5J^3-\tfrac{48}5J|V|^2
\right)d\x\right)_\w-108.\label{expocls2}
\end{gather}
Each may be estimated to be
nonnegative, with equality if\/f $e^\w$ is a conformal factor.  Thus
we ``win'', for a given operator $A$,
if there are no conf\/licting signs
in the expression of the $\cD$-functional for $A$ as
\[
a_0\cF_0+a_1\cF_1+a_2\cF_2+a_3\cF_3.
\]
Remarkably, there are no conf\/licting signs in the cases of
$Y$ and $\dc^2$:
up to constant multiples, $(a_0,a_1,a_2,a_3)$ is{\samepage
\begin{alignat*}{3}
&Y:&& (600,6,23,10),& \\
&\dc^2: \ \ && -(11460,93,556,365).&
\end{alignat*}

}

The nonnegativity of $\cF_1(\w)$ is clear; it vanishes if\/f $J$ is constant,
if\/f (by Obata's Theorem) $\w$~is a conformal factor.  An analysis
of $\cF_2(\w)$ may be done by considering the Yamabe quotient applied
to the function $J$.  The Yamabe quotient is
\[
\cY(u)=\dfrac{\int u(\D+2J)u}{\|u\|_3^2},
\]
so the Yamabe constant (minimum of the Yamabe quotient) on
${\bf S}^6$ is 6 (setting $J=3$ and $u=1$).
This gives
\[
\int_{{\bf S}^6}\left(\left\{|dJ|^2+2J^3\right\}\right)_\w
\ge 6\left(\int_{{\bf S}^6}(|J|^3d\x)_\w)\right)^{2/3},
\]
since the Yamabe constant is conformally invariant (so in particular,
it is the same at $g_\w$ as at~$g_0$).
Now run an argument analogous to \nnn{hold}, based on
\nnn{beclp} below, and replacing the
Schwarz inequality step with an application of the H\"older inequality,
with exponents~$\frac32$ and~3.
(A very general version of this argument is given in \cite[Corollary 3.6]{sharp}.)
The result is the nonnegativity of~$\cF_2(\w)$, with equality if\/f
we are at a Yamabe metric, if\/f (by Obata's Theorem) $\w$ is a conformal
factor.

The estimate of $\cF_3(\w)$ is more subtle.  By \cite{wunsch86}, there
is a second-order conformal covariant~$D$ on trace-free
symmetric 2-tensors.  One may compute \cite{sharp} that
for $n\ge 5$, this operator
has positive def\/inite leading symbol, is formally self-adjoint, and has
has positive spectrum on the round sphere.
Since the sign of the bottom eigenvalue (if any) of a~conformal covariant
is a~conformally invariant quantity, $D$~has positive spectrum
for metrics
conformal to the round one.  Thus for $n\ge 5$,
\[
\int_{{\bf S}^n}(D_\w\f,\f)(dv)_\w\ge 0,
\]
with equality if\/f $\f=0$.  If we substitute the trace-free Ricci tensor
$b_\w$
of the metric $g_\w$
for $\f$, the result (when $n=6$) is the inequality $\cF_3(\w)\ge 0$.
Furthermore, we have equality if\/f $b_\w=0$, if\/f $({\bf S}^6,g_\w)$ is a space
of constant curvature, if\/f $\w$ is a conformal factor.

An interesting aspect of the determinant calculation for the conformal
Laplacian in dimension 6 is the appearance of certain $n-8$ factors.
In fact,
\begin{gather*}
(4\p)^{n/2}7!a_6[Y]=(n-8)\int\left\{
-3(n-6)\njs
-\tfrac19(35n^2-266n+456)J^3\right. \\
\left.\phantom{(4\p)^{n/2}7!a_6[Y]=}{} +\tfrac23(n-1)(7n-30)J|V|^2
-\tfrac29(5n^2-2n-48)\tr V^3\right\}dv \\
\phantom{(4\p)^{n/2}7!a_6[Y]}{}=-\tfrac53(n-8)\int Q_6dv
+(n-8)(n-6)\int\left\{
-\tfrac{13}6\njs\right. \\
\left.\phantom{(4\p)^{n/2}7!a_6[Y]=}{}-\tfrac1{36}(125n-314)J^3
+\tfrac23(7n-5)J|V|^2
-\tfrac29(5n+28)\tr V^3\right\}dv.
\end{gather*}

This is analogous to an $n-6$ factor that appeared in the 4-dimensional
case.  The root cause of this phenomenon
is the relative conformal invariance
of the heat invariant $U_{n-2}[Y]$.
To our knowledge, this was f\/irst remarked
upon by Schimming \cite{schim}, in connection with studies of
Huygens' principle for hyperbolic equations.  A Riemannian geometry
proof is
given in \cite{pr87}, based on the fact that
$U_{n-2}[Y]$ is the coef\/f\/icient of the ``f\/irst log term'' in the
radial expansion of the fundamental solution (Green's function) of $Y$.
\cite{bgcpde} shows that this phenomenon is really a~consequence of
the conformal index property: take a f\/ixed operator of the form
\mbox{$D=-\nd^a\nd_a-E$}, and consider perturbations
\[
e^{-2\e\w}(D-\d F),
\]
where $\e$ and $\d$ are real parameters and $\w$ and $F$ are functions.
Now compute
\[
\left.\dfrac{\partial^2}{\partial\d\partial\e}\right|_{(\d,\e)=(0,0)}
\]
of the heat expansion
with each of the two possible partial derivative orderings; the result
is the relative conformal invariance of $U_{n-2}[Y]$:
\[
U_{n-2}[Y]_\w=e^{-(n-2)\w}U_{n-2}[Y]_0\,.
\]
(That is, the conformal variation of $e^{(n-2)\w}U_{n-2}[Y]_\w$
vanishes.)
In the conformally f\/lat category, this forces $U_{n-2}[Y]$ to vanish.
Thus, for example, $U_6[Y]$ is forced to have an overall factor
of $n-8$.

This calculation gives no special conclusion for powers of a conformal
covariant (like $\dc^2$).  Replacing $Y$ by
a higher-order operator
leads to the relative conformal invariance of a dif\/ferent term
in the heat asymptotics; for example, $U_{n-4}[P]$ in the case of
the Paneitz operator $P$.

Recently, Larry Peterson and the author \cite{bp} have
computed the determinant quotient for the conformal Laplacian in
the standard conformal class on ${\bf S}^8$.  The starting point is a
general formula for $U_8[D]$ due to Avramidi \cite{avra}, where
$D$ is any operator of Laplace type.  We use this to compute
$a_8[Y]$ in arbitrary dimension $n$, in the conformally f\/lat category.
The result is
$(4\p)^{-n/2}$ times a linear combination (with rational-in-$n$ coef\/f\/icients)
of the integrals of the nine scalar invariants
\[
(\D J)^2,\quad J|dJ|^2,\quad J|\N V|^2,\quad (V,dJ\otimes dJ),\quad J^4,\quad
J^2|V|^2,\quad J{\rm tr}\,V^3,\quad |V|^4,\quad {\rm tr}\,V^4.
\]
Here, for example, ${\rm tr}\,V^3=V^a{}_bV^b{}_cV^c{}_a$.
Though this is a long calculation, something of a check is performed
by the expected $n-10$ factor -- the fact that this factor appears is
reassuring.  We have not yet been able to use this to solve the
extremal problem for $\cD(Y,g_\w)$, but several sharp inequalities
are available.

An interesting aspect of the extremal problem on low-dimensional
spheres is the ``checkerboard pattern''
\def\skemastrut{\vrule height 0.4 cm depth 0.2 cm width 0 cm}
\def\skema#1\endskema{\vbox{\def\\{\skemastrut\cr\noalign{\hrule}}%
   \offinterlineskip
   \halign{\vrule\hskip1em
  ##\hskip1em\hfil\vrule&&\hfil\quad$##$\quad\hfil\vrule\cr
\noalign{\hrule}
  #1\skemastrut\crcr\noalign{\hrule}}}}
\vskip -0.4cm
$$
\skema
 & \det Y & \det\dc^2 \\
${\bf S}^2$ & \max & \min \\
${\bf S}^4$ & \min & \max \\
${\bf S}^6$ & \max & \min
\endskema
$$
Here ``max'' means that the (suitably penalized) determinant quantity
attains a max at the
at the round metric and its conformal dif\/feomorphs.

In (as yet inconclusive) work with Carlo Morpurgo and Bent \O rsted,
we have attempted to explain at least the max/min alternation for
the conformal Laplacian as follows.  Suppose $n$ is even.
Morpurgo shows that
for real $s$ in the interval $(\frac{n}2-1,\frac{n}2)$, the zeta
dif\/ference
\begin{equation}\label{zdif}
\z_{Y_\w}(s)-\z_{Y_0}(s)
\end{equation}
is nonnegative.  This dif\/ference has potential simple poles
at the positive integers $\le n/2$, so if all these poles
are realized, and no zeros intervene, we have a quantity of
sign $(-1)^{\frac{n}2-1}$ for small positive $s$; by the conformal
index property, this should also be the sign of the log-determinant
dif\/ference
\begin{equation}\label{lddif}
\z'_{Y_\w}(0)-\z'_{Y_0}(0).
\end{equation}
Recalling the minus sign in the def\/inition of the determinant as
$\exp(-\z'(0))$, this agrees with the $n=2,4,6$ cases
in the above table.

Considering what happens as we move along curves in the
conformal class, by the obvious ref\/lection
principle, intervening zeros of \nnn{zdif} on the real axis must
either come in pairs (i.e.\ be double zeros), or must move along the
real axis.  With double zeros, no harm is done to the above argument.
A single zero, though, has the potential to reverse the sign of the
above result, or ruin the result entirely, if it moves past $s=0$.
As it moves through $s=0$, what we see is a zero of the
log-determinant dif\/ference \nnn{lddif}.  So the problem is to rule
these out -- clearly a~restatement of the assertion that the
round metric is an absolute extremal.  (Among the many things needed
to make this discussion precise is to account for the conformal
dif\/feomorphisms, which preserve the zeta function.)

In any case, what is missing is some sort of general-dimension approach
to the convexity properties of the determinant quotient.
This would avoid the dimension-by-dimension approach, and have the
side benef\/it of {\em proving} sharp harmonic-analytic inequalities
rather than just using them.  This is conceptually related to
the topic of {\em complementary series} in representation theory,
as we shall try to explain in the next section.

\section{The complementary series}
The constants and dif\/ferential operators that appear
in the sharp inequalities
that estimate quantities like those in \nnn{expocls}
and \nnn{expocls2}, and in turn the determinant, also appear
in the study of the {\em complementary series} of
SO${}_0(n+1,1)$.

The conformal transformations $h$ of any Riemannian manifold,
{\rm ctran}$(M,g)$, form a group.  Recall that
such a transformation has
\[
h\cdot g=(h^{-1})^*g=\W_h^2g,\qquad 0<\W_h\in\ci(M).
\]
The corresponding inf\/initesimal notion is that of {\em conformal vector
fields} $X$; these satisfy
\[
\cL^{}_Xg=2\w^{}_Xg,\qquad \w^{}_X\in\ci(M),
\]
where $\cL$ is the Lie derivative,
and form a Lie algebra ${\rm cvf}(M,g)$.
The {\em conformal factors} satisfy {\em cocycle conditions}
\[
\W_{h_1\circ h_2}=\W_{h_1}(h_1\cdot\W_{h_2}), \qquad
\w_{[X_1\,,X_2]}=X_1\w_{X_2}-X_2\w_{X_1}.
\]
The cocycle conditions are equivalent to the assertion that the
family ({\em series}, in representation theoretic parlance) of maps
\[
u_a(h)=\W_h^ah\cdot, \qquad
U_a(X)=\cL_X+a\w_X
\qquad (a\in\bbC)
\]
are homomorphisms into the group of automorphisms (resp.\ Lie algebra
of endomorphisms) of functions on $M$, or in fact
tensors, or tensor-spinors, of any type.
The {\em isometry} subgroup (subalgebra) is def\/ined by the condition
$\W_h=1$ ($\w_X=0$).

When $(M,g)$ is standard ${\bf S}^n$, these
objects give the principal series of SO${}_0(n+1,1)$ (or Spin${}_0(n+1,1)$
if spinors are involved).  That is,
\[
{\bf S}^n=G/MAN=K/M,
\]
where
\[
G={\rm SO}_0(n+1,1),\quad M={\rm SO}(n),\quad A\cong\bbR_+,
\quad N\cong\bbR^n,\quad
K={\rm SO}(n+1).
\]

The representations are ${\rm Ind}_{MAN}^G\l\otimes\nu\otimes 1$,
with the interpretations
\begin{gather*}
\l\approx{\rm tensor}\!-\!{\rm spinor\ species}; \qquad
\nu\approx{\rm conformal\ weight}.
\end{gather*}
The stereographic map $\bbR^n\hookrightarrow
\bbS^n$ is a special case of the more
general Iwasawa decomposition theoretic map
\[
\bar N\hookrightarrow G=\bar NMAN\to
G/MAN.
\]

Such induced representations
are best thought of {\em algebraically} at f\/irst,
as $({\mathfrak{g}},K)$-modules.  This means: take the $K$-f\/inite vectors,
so one has a direct sum of $K$-modules.  The action of $\mathfrak{g}$
lets you {\em travel} among these modules -- you can get from
a $K$-module $\a$ to a $K$-module $\b$ only if $\mathfrak{s}\otimes\a$
has a $\b$ summand in its $K$-decomposition (the {\em selection rule}).
Here
\[
\mathfrak{g}=\mathfrak{k}+\mathfrak{s}
\]
is the Cartan decomposition.

An example is the expansion into spherical harmonics of functions
on $\bbS^n$:
\[
\cE\cong_{{\rm SO}(n)}E_j.
\]
Using the {\em big} Lie algebra ${\mathfrak{so}}(n+1,1)$, we can move
up and down the ladder,
\[
E_{j-1}\leftarrow E_j\rightarrow E_{j+1},
\]
within the {\em spherical principal series} representations
\[
{\rm Ind}_{MAN}^G 1\otimes\nu\otimes 1=:(U_{\frac{n}2+\nu},u).
\]
The unitaries in this series are at imaginary $\nu$
(the {\em unitary principal series}), and, on the real axis,
in the interval $|\nu|<n/2$ -- the {\em complementary series}.
The inner products on the complementary series representations
$(U_{\frac{n}2+\nu},u)$
have the form
\[
(\f,\y)_\nu=\int_{\bbS^n}\f(A_{2\nu}\bar\y)d\x,
\]
where $A_{2\nu}$ is the {\em intertwining operator}
\[
A_{2\nu}:(U_{\frac{n}2+\nu}\,,u)\to(U_{\frac{n}2-\nu},u).
\]
Since the two representations $U_{\frac{n}2\pm\nu}$
are naturally dual (they live in section spaces of the bundles
of $\mp\nu$-densities), the integrand is natural.

The Knapp--Stein intertwinor \cite{ks80} is an integral operator,
the integral converging in a certain range of Re$(\nu)$.
Dif\/ferential intertwinors, like the realization
of the conformal Laplacian
and Paneitz operator on the sphere, live in the analytic continuation
(in $\nu$) of the family of Knapp-Stein intertwinors.
In the sense of pseudo-dif\/ferential operators and of complex powers
of an elliptic operator \cite{seeley}, assuming this analytic continuation,
\[
A_{2\nu}=\D^\nu+{\rm(lower\ order)}.
\]
Up to a constant, $(P_m)_0$ is $A_{m}$.
Not surprisingly, the obvious continuation of the
formula \nnn{old87}
given earlier for $(P_m)_0$ turns out to be a formula
for $A_{2\nu}$ (all $\nu$):
\[
A_{2\nu}=\dfrac{\G(\cB+\nu+\frac12)}{\G(\cB-\nu+\frac12)},
\]
where (recall)
\[
\cB=\sqrt{\D+\left(\dfrac{n-1}2\right)^2}.
\]
For $|\nu|<n/2$, by Stirling's formula for the Gamma function,
the complementary series norm
on $U_{\frac{n}2+\nu}$ is
a Sobolev $L^2_{\nu}$ norm.

The fact that the complementary series norm is an {\em invariant}
Sobolev norm follows from the intertwining property:
\begin{gather*}
\dint_{\bbS^n}\f\cdot A_{2\nu}\left\{\left(\cL_X+\left(\frac{n}2-\nu\right)
\w_X\right)
\y\right\}d\x
=\dint_{\bbS^n}\f\cdot
\left(\cL_X+\left(\frac{n}2+\nu\right)\w_X\right)
A_{2\nu}\y\,d\x \\
\qquad{}
=\dint_{\bbS^n}\left\{\left(-\cL_X-\left(\frac{n}2-\nu\right)\w_X\right)\f
\right\}\cdot
A_{2\nu}\y d\x.
\end{gather*}
Here we're using the fact that the case $\nu=0$ is in the
unitary principal series; this really follows from the simpler fact
that the volume form $E$ has
\[
\cL_XE=n\w^{}_XE.
\]
This last formula may be viewed as the origin of the rho-shift $n/2$.
This also happens more generally in the theory of parabolic induction for
semisimple Lie groups; the rho-shift due to the $({\mathfrak{g}},{\mathfrak{a}})$
root system is related to the response of the natural volume form
on $G/MAN$ to inf\/initesimal transformations from ${\mathfrak{g}}$.
In the present case, the calculation is
\begin{gather*}
0=\dint\underbrace{\cL_X}_{d\iota(X)}\left(\f\cdot A_{2\nu}\y\,d\x\right)
=\dint\f\Bigg\{\underbrace{\left(\cL_X+\left(\frac{n}2+\nu\right)
\w_X\right)A_{2\nu}}_{A_{2\nu}\left(\cL_X+\left(\frac{n}2-\nu\right)
\w_X\right)}\y
\Bigg\}d\x \\
\phantom{0=}{}
+\dint\left\{\left(\cL_X-\left(\frac{n}2+\nu\right)\w_X\right)\f\right\}
A_{2\nu}\y\cdot d\x
+\dint\f(A_{2\nu}\y)n\w_X\,d\x \\
\phantom{0}=\dint\f\Bigg\{
A_{2\nu}\left(\cL_X+\left(\frac{n}2-\nu\right)
\w_X\right)\y
\Bigg\}d\x
+\dint\left\{\left(\cL_X+\left(\frac{n}2-\nu\right)\w_X\right)\f\right\}
A_{2\nu}\y\cdot d\x.
\end{gather*}

There's another natural invariant norm on the space carrying
$U_{\frac{n}2-\nu}$,
namely the
$L^p$ norm with
\[
p=\dfrac{2n}{n+2\nu}.
\]
Note that
\[
p\in(1,\infty)\iff\nu\in\left(-\dfrac{n}2\,,\dfrac{n}2\right).
\]
The invariance calculation is
\begin{gather*}
0=\dint\cL_X\left(\f^pd\x\right)
=p\dint\f^{p-1}(\cL_X\f)d\x+\dint\f^pn\w_Xd\x
=p\dint\f^{p-1}\left(\cL_X+\dfrac{n}{p}\w_X\right)f\,d\x;
\end{gather*}
this shows that when
\[
\dfrac{n}{p}=\dfrac{n}2+\nu,
\]
we have an invariant norm:
\[
\dint\f^p\,d\x=\dint h\cdot(\f^p\,d\x)=\dint(h\cdot\f)^p\,
\underbrace{h\cdot d\x}_{\W_h^n\,d\x}
=\dint\W_h^{-p\left(\frac{n}2+\nu\right)}(u_\nu(h)\f)^p\W_h^n\,d\x.
\]

The Sobolev embedding inequalities
\[
L^2_\nu\hookrightarrow L^{\frac{2n}{n-2\nu}},\qquad
L^{\frac{2n}{n+2\nu}}\hookrightarrow L^2_{-\nu}
\qquad (\nu\ge 0),
\]
in the sharp form due to Beckner \cite{bec}, compare these two
invariant norms for $\nu$ between $-n/2$ and 0, saying that
saying that
\begin{equation}\label{beclp}
\displaystyle\max_{L^p\setminus 0}\dfrac{|(\f,\y)_\nu|}{\|\f\|_p\|\y\|_p}
\end{equation}
is attained exactly at conformal factors: both $\f$ and $\y$
should be constant multiples of $\W_h^{\frac{n}2-\nu}$, where
$h$ is a conformal transformation with $(h^{-1})^*g=\W_h^2g$.
This gives a comparison of norms, and of norms on the corresponding
dual spaces.
These lead, by endpoint dif\/ferentiation, to the exponential class
inequality \nnn{expocls}.
They are also used in deriving geometric inequalities like
the estimate on $\cF_2(\w)$ in \nnn{expocls2}, based on
\cite[ Corollary 3.6]{sharp}.

This interpretation in terms of representation theory is very
clean and nice, but unfortunately, it has never been put to work
in actually {\em proving} sharp inequalities.  This remains a
prospect for the future, much like (and probably closely related to)
that of proving convexity properties of the determinant,
and using them to solve extremal problems.  The $L^p$ norms
in \nnn{beclp} are quite amenable to hard geometro-analytic
methods like symmetric decreasing rearrangement, while
representation theory is still quite {\em linear} (meaning
{\em quadratic} in the present formulation -- that is, concerned
with inner products).

\pdfbookmark[1]{References}{ref}
\LastPageEnding

\end{document}